\documentclass[10pt,a4paper]{article}%
\usepackage[centertags]{amsmath}
\usepackage{amsfonts}
\usepackage{amssymb}
\usepackage{amsthm}
\usepackage{epsfig}
\usepackage{setspace}
\usepackage{ae}
\usepackage{eucal}
\usepackage[usenames]{color}%
\usepackage{graphicx}

\theoremstyle{plain}
\newtheorem{thm}{Theorem}[section]
\newtheorem{lem}[thm]{Lemma}
\newtheorem{cor}[thm]{Corollary}
\newtheorem{prop}[thm]{Proposition}
\theoremstyle{definition}
\newtheorem{defi}[thm]{Definition}

\theoremstyle{remark}

\newtheorem{rem}[thm]{Remark}

\begin{document}

\title{On global regularity and singularities of  Navier-Stokes- and Euler equation solutions }
\author{J\"org Kampen }
\maketitle

\begin{abstract}
Euler-Leray data functions of first and second order are defined by first and second order derivatives of the nonlinear spatial part of the incompressible Euler equation operator in Leray projection form applied to Cauchy data. The Lipschitz continuity (or even a strong H\"{o}lder continuity) of these functions for strong Cauchy data in $H^m\cap C^m,~m\geq 2$ is sufficient for the existence of global regular upper bounds of incompressible Navier Stokes equation solutions (in case of Cauchy data in $H^m\cap C^m,~m\geq 2$). Global regular upper bounds of global solution branches of the incompressible Euler equation can be obtained (in case of Cauchy data in $H^m\cap C^m,~m\geq 3$)), if the Cauchy data satisfy an additional condition of strong polynomial decay at spatial infinity. Furthermore, if a Lipschitz condition for the Euler-Leray data function of second order is satisfied, then there are long time vorticity blow ups of the incompressible Euler equation, and, correspondingly, short time and long time vorticity blow ups or singular solutions of incompressible Navier Stokes equations with time-dependent forces in $H^{m-2}\cap C^{m-2}$. A further consequence is that multiple global solutions of the incompressible Euler equations exist. The so-called degeneracy issue of convolutions with first order spatial derivatives of the Gaussian is addressed. For the Navier Stokes equation Lipschitz continuity of the Leray projection term in local time has the effect that a rough upper bound estimate of the Leray projection term in terms of the second moment of the Gaussian implies the existence of global regular upper bounds due to stronger damping of spatially scaled solutions. These simple estimates can be refined in order to obtain global regular upper bounds of solution branches in the viscosity limit. 
%Global regularity and long time singularities of Navier-Stokes- and Euler equation solutions
\end{abstract}

%\footnotetext[1]{Weierstrass Institute for Applied Analysis and Stochastics,
%M%ohrenstr.39, 10117 Berlin, Germany. Support by DFG Matheon is acknowledged.
%\texttt{{kampen@wias-berlin.de}}.}

\section{Euler Leray data functions of first order and global regular upper bounds of the Navier Stokes equation}

Consider the Cauchy problem for the incompressible Navier Stokes equation for regular velocity data $v^{\nu}_i(t_0,.)\in H^m\cap C^m,~m\geq 2,~1\leq i\leq D$ on a domain $\left[t_0,t_0+\Delta\right]\times {\mathbb R}^D$  at time $t_0\geq 0$ and for a small time horizon $\Delta>0$. The short time solution of the incompressible Navier Stokes equation Cauchy problem with constant positive viscosity $\nu >0$ and without external forces has the representation
\begin{equation}\label{Navlerayscheme}
\begin{array}{ll}
 v^{\nu}_i=v^{\nu}_i(t_0,.)\ast_{sp}G_{\nu}
-\sum_{j=1}^D \left( v^{\nu}_j\frac{\partial v^{\nu}_i}{\partial x_j}\right) \ast G_{\nu}\\
\\+\left( \sum_{j,m=1}^D\int_{{\mathbb R}^D}\left( \frac{\partial}{\partial x_i}K_D(.-y)\right) \sum_{j,m=1}^D\left( \frac{\partial v^{\nu}_m}{\partial x_j}\frac{\partial v^{\nu}_j}{\partial x_m}\right) (.,y)dy\right) \ast G_{\nu}.
\end{array}
\end{equation}
Here, the symbol '$\ast$' denotes convolution with respect to space and time and the symbol '$\ast_{sp}$' denotes convolution with respect to the spatial variables. Furthermore, the symbol $K_{D}$ refers to the the Laplacian kernel of dimension $D$, and the symbol $G_{\nu}$ denotes the Gaussian fundamental solution of the heat equation
\begin{equation}\label{fund1}
p_{,t}-\nu \Delta p=0.
\end{equation} 
The evaluation of such a fundamental solution is usually denoted by $G_{\nu}(t,x;s,y)$ where in the special $s=0$ and $y=0$ we also write
\begin{equation}
G_{\nu}(t,x):=G_{\nu}(t,x;0,0) \mbox{ for the sake of brevity.}
\end{equation}
The validity  of the local-time representation in (\ref{Navlerayscheme}) is due to a local time contraction result with respect to the norm $\max_{1\leq i\leq D}\sup_{t\in [t_0,t_0+\Delta]}{\big |}v^{\nu}_{i}(t,.){\big |}_{H^m\cap C^m}$ (for a small time interval size $\Delta_0 >0$ and regularity order $m\geq 2$). The solution has then an upper bound with respect to this norm on the time interval $[t_0,t_0+\Delta_0]$.
We may use the incompressibility condition
\begin{equation}\label{incomp}
\sum_{j=1}^D\frac{\partial v^{\nu}_j}{\partial x_j}=0,
\end{equation}
in order to rewrite the Burgers term. Indeed, the incompressibility condition in (\ref{incomp}) implies that
\begin{equation}
\sum_{j=1}^D\frac{\partial ( v^{\nu}_iv^{\nu}_j)}{\partial x_j}=\sum_{j=1}^Dv^{\nu}_j\frac{\partial v_i}{\partial x_j}+v_i\sum_{j=1}^D\frac{\partial v^{\nu}_j}{\partial x_j}= \sum_{j=1}^Dv^{\nu}_j\frac{\partial v^{\nu}_i}{\partial x_j}.
\end{equation}
Hence we may rewrite (\ref{Navlerayscheme}) such that the nonlinear terms are convolutions with first order spatial derivatives of the Gaussian. We have
\begin{equation}\label{Navlerayscheme2i}
\begin{array}{ll}
 v^{\nu}_i=v^{\nu}_i(t_0,.)\ast_{sp}G_{\nu}
-\sum_{j=1}^D \left( v^{\nu}_jv^{\nu}_i\right) \ast G_{\nu,j}\\
\\+\left( \int_{{\mathbb R}^D}\left( K_D(.-y)\right) \sum_{j,m=1}^D\left( \frac{\partial v^{\nu}_m}{\partial x_j}\frac{\partial v^{\nu}_j}{\partial x_m}\right) (.,y)dy\right) \ast G_{\nu,i},
\end{array}
\end{equation}
where we use the convolution rule for derivatives.  We shall observe that the viscosity damping of the term $v^{\nu}_i(t_0,.)\ast_{sp}G_{\nu}$ in 
(\ref{Navlerayscheme2i}) is stronger than possible growth of the nonlinear terms due to the spatial effects related to $G_{\nu,j}$ if the data ${\big |}v^{\nu}_i(t_0,.){\big |}_{H^m\cap C^m},~m\geq 2$ exceed a certain level, and if the Leray projection term satisfies a Lipschitz condition. For the scaled Gaussian $G^{\rho,r}_{\nu}$ 
(obtained by replacement of $\nu$ by $\rho r^2\nu$, cf. below), we compute
\begin{equation}
\begin{array}{ll}
{\big |}G^{\rho,r}_{\nu,i}(\tau,y){\big |} 
={\Big |}\frac{-2y_i}{4 \rho r^2\nu \tau}\frac{1}{\sqrt{4\pi\rho r^2\nu t}^D}\exp\left(-\frac{|y|^2}{4\rho r^2\nu \tau} \right){\Big |}\\
\\
\leq \frac{2}{\sqrt{\pi}^D}\frac{1}{(4\nu\rho r^2 \tau)^{\delta}|y|}\left( |y|^2\right)^{\delta-D/2} \left( \frac{|y|^2}{4\rho r^2\nu \tau}\right)^{D/2+1-\delta} \exp\left(-\frac{|y|^2}{4\rho r^2\nu \tau} \right).
\end{array}
\end{equation}
Hence, we have for $\delta \in (0,1)$ and all $\rho,r>0$ 
\begin{equation}\label{grrh0}
{\big |}G^{\rho,r}_{\nu ,i}(\tau,y){\big |}\leq \frac{C}{(4 \rho r^2 \nu \tau)^{\delta}|y|^{D+1-2\delta}}, 
\end{equation}
where the upper bound constant $$C=\sup_{|z|>0}\left( z\right)^{D/2+1-\delta} \exp\left(-z^2\right) >0$$ is sufficient and independent of $\nu >0$. The estimate in (\ref{grrh0}) may be used locally, i.e., on compact subsets of ${\mathbb R}^D$. The symmetry of spatial first order derivatives of the Gaussian, i.e., the relation
\begin{equation}
G_{\nu,i}(t,y)=\frac{-2y_i}{4 \rho r^2\nu \tau}\frac{1}{\sqrt{4\pi\rho r^2\nu t}^D}\exp\left(-\frac{|y|^2}{4\rho r^2\nu \tau} \right)=-G_{\nu,i}(t,y^{i,-})
\end{equation}
for $y^{i,-}=\left(y^{i,-}_1,\cdots,y^{i,-}_n\right)$ with $y^{i,-}_j=y_j$ for $j\neq i$ and $y^{i,-}_i=y_i$ can be exploited in convolutions with Lipschitz continuous functions or even in convolutions with H\"{o}lder continuous functions.
Concerning the Leray projection data term in (\ref{Navlerayscheme2i}) and in similar representations for multivariate spatial derivatives up to order $m$ we shall observe that Lipschitz continuous or H\"{o}lder continuous Leray data functions are convoluted (with respect to space and time) with first order spatial derivatives of the Gaussian.

For a Lipschitz continuous function $L_x=L(\tau,x-y)$ with $|L_x(\tau,y)-L_x(\tau,y')|\leq l|y-y'|$ (for a Lipschitz constant $l$ which is independent of $x$) and on a small time interval of length $\Delta_0$ we first (for any given $x\in {\mathbb R}^D$) get the local estimate
\begin{equation}\label{lipob}
\begin{array}{ll}
\int_{t_0}^{t_0+\Delta_0}{\big |}\int_{|y|\leq \sqrt{4\rho r^2\nu} } L(\tau,x-y)
G^{\rho,r}_{\nu,i}(t_0+\Delta_0-\tau,y){\big |}dyd\tau\\
\\
=\int_{t_0}^{t_0+\Delta_0}{\Big |}\int_{|y|\leq \sqrt{4\rho r^2\nu} }{\Big |}L_x(\tau,y)-L_x(\tau,y^{i,-}){\Big |}G^{\rho,r}_{\nu,i}(t_0+\Delta_0-\tau,y)dy{\Big |}d\tau \\
\\
\leq \int_{t_0}^{t_0+\Delta_0}{\Big |}\int_{|y|\leq \sqrt{4\rho r^2\nu} }{\Big [} \frac{4|y|^2}{4 \rho r^2\nu (t_0+\Delta_0-\tau)}\exp\left(-\frac{|y|^2}{8\rho r^2\nu (t_0+\Delta_0-\tau)} \right)\\
\\
\frac{1}{\sqrt{4\pi\rho r^2\nu(t_0+\Delta_0- \tau)}^D}\exp\left(-\frac{|y|^2}{8\rho r^2\nu (t_0+\Delta_0-\tau)} \right){\Big ]}dy d\tau\\
\\
\leq \int_{t_0}^{t_0+\Delta_0}{\Big |}\int_{|y|\leq \sqrt{4\rho r^2\nu} }l\frac{C}{(4 \rho r^2 \nu (t_0+\Delta_0-\tau)^{\delta}|y|^{D-2\delta}}{\Big |}dy d\tau\\
\\
\leq lC\Delta_0^{1-\delta}\int_{|y|\leq \sqrt{4\rho r^2\nu}} \frac{|y|^{2\delta}}{4 \rho r^2 \nu}dy\\
\\
\leq lC\Delta_0^{1-\delta}\frac{|y|^{2\delta}}{4 \rho r^2 \nu}{\Big |}^{\sqrt{4\rho r^2\nu}}_0\leq lC^*\Delta_0^{1-\delta}
\end{array}
\end{equation}
for some finite constants $C$ and $C^*$. Furthermore the convoluted Gausssian 
\begin{equation}
\frac{1}{\sqrt{4\pi\rho r^2\nu(t_0+\Delta_0- \tau)}^D}\exp\left(-\frac{|y|^2}{8\rho r^2\nu (t_0+\Delta_0-\tau)} \right)~ 
\end{equation}
is of polynomial decay of any order at spatial infinity. Especially for a $L(.-y)\in H^m$ we have
\begin{equation}
\int_{t_0}^{t_0+\Delta_0}{\big |}\int_{|y|\leq \sqrt{4\rho r^2\nu} } L(\tau,.-y)
G^{\rho,r}_{\nu,i}(t_0+\Delta_0-\tau,y){\big |}dyd\tau\in H^m\cap C^2.
\end{equation}
 
We shall refine these estimates below and show how they work in the context of $H^m\cap C^m$-norms. 
Note that the constants $C,C^*>0$ above can be chosen independent of the parameters and the viscosity as we integrate spatially up to $\sqrt{4\rho r^2\nu}$.
For the scaled Navier Stokes equation, i.e. the Navier Stokes equation in coordinates $(\tau,y)$ with $t=\rho \tau$ and $x=r y$  the nonlinear get a factor $\rho r$ such hat we get an upper bound proportional to $\rho r\Delta_0^{1-\delta}$, while we shall observe a damping of order $\nu\rho r^2\Delta_0^3$ on a small time scale and for data which exceed a certain threshold. This implies that upper bounds can be constructed that are preserved by a scheme for some strong spatial scaling $r>1$ if $\nu$ is positive.
Indeed we shall observe below that this construction idea or a refinement of this construction idea is sufficient in order to get global regular upper bounds for the velocity component solution functions $v_i(t,.)=v^{\rho,r}(\tau,.)$ in the function space $H^m\cap C^m$ for $m\geq 2$ which are preserved from time $t_0$ to time $t_0+\Delta$, where the time interval length $\Delta$ with respect to original time coordinates $t$ is related to the time interval length $\Delta=\rho \Delta_0$ to the time interval length $\Delta_0$ in scaled coordinates. An extension of the argument (with some restrictions concerning regularity) to the incompressible Euler equation is obtained by compactness arguments and upper bound constructions for vorticity equations. This leads to global classical solution branches for the Euler equations. In the case of an incompressible Euler equation upper bounds for stronger norms can be obtained by auto-controlled schemes (cf. below). Note that the construction is not straightforward because for $y=rx$ the derivative of $v^{r}_i(t,y):=v_i(t,x)$ becomes
\begin{equation}
v_{i,j}(t,x)=r v^{r}_{i,j}(t,y).
\end{equation}
This means that for a strong scaling $r>1$ we have 
\begin{equation}\label{normcomp}
{\big |} v_{i,j}(t,.){\big |}_{H^m\cap C^m}\leq c(n)r^m{\big |} v^r_{i,j}(t,.){\big |}_{H^m\cap C^m}.
\end{equation}
Hence the spatial scaling parameter $r$ should be independent of $\nu$ in a construction of an uniform regular upper bound for a solution of the incompressible Euler equation. We do not have this restriction for the auto-controlled schemes where upper bounds are constructed which are linear in time.  
First we we shall obtain upper bounds for the right side of (\ref{normcomp}) as $r>0$ becomes large and for positive $\nu >0$. Although the upper bound consruction is independent of $\nu >0$ for the value function $v^{r}_i,~1\leq i\leq n$ itself, it is not possible to construct regular upper bounds which hold in the viscosity limit this way. This is clear as a reciprocal dependence of $r$ and (some positive power of $\nu$) leads to exploding data. However it is possible to apply our technique to certain systems of spatial derivatives of the value function $v^{r}_{i},~1\leq i\leq n$ if there is a similar spatial  scaling of derivative Burgers and Leray projection terms compared to the viscosity (Laplacian) term as for the original equation. For example this is true for the vorticity equation. 
In coordinates the vorticity Euler equation is the viscosity limit of a family of equations of the form 
\begin{equation}\label{vort}
 \frac{\partial \omega^{\nu}_i}{\partial t}-\nu \Delta \omega^{\nu}_i+\sum_{j=1}^3v^{\nu}_j\frac{\partial \omega^{\nu}_i}{\partial x_j}=\sum_{j=1}^3\frac{1}{2}\left(\frac{\partial v^{\nu}_i}{\partial x_j}+\frac{\partial v^{\nu}_j}{\partial x_i} \right)\omega^{\nu}_j.
\end{equation}
Here $v^{\nu}_i,¸1\leq i\leq D$ denotes the original velocity related to a positive viscosity $\nu >0$, where we have
\begin{equation}\label{vel}
v^{\nu}(t,x)=\int_{{\mathbb R}^3}K_3(x-y)\omega^{\nu}(t,y)dy,~\mbox{where}~K_3(x)h=\frac{1}{4\pi}\frac{x\times h}{|x|^3}.
\end{equation}
Note that upper bounds for first order spatial derivatives $v^{\nu}_{i,j},~1\leq i,j\leq D$ can be obtained from upper bounds of the vorticity $\omega^{\nu}_{i},~1\leq i\leq D$. We shall use data of strong spatial polynomial decay in order to transformation to compact spaces (cf. below). This is useful whether a Rellich compactness argument or convergence in another Banach space is considered. Since there is a small loss of regularity in the situation of Rellich's theorem (which we can still use in a variation of argument where we use auto-controlled schemes), for a direct argument we shall consider a classical Banach space for compactness arguments. In this context recall
\begin{prop}
For open and bounded $\Omega\subset {\mathbb R}^n$ and consider the function space
\begin{equation}
\begin{array}{ll}
C^m\left(\Omega\right):={\Big \{} f:\Omega \rightarrow {\mathbb R}|~\partial^{\alpha}f \mbox{ exists~for~}~|\alpha|\leq m\\
\\
\mbox{ and }\partial^{\alpha}f \mbox{ has an continuous extension to } \overline{\Omega}{\Big \}}
\end{array}
\end{equation}
where $\alpha=(\alpha_1,\cdots ,\alpha_n)$ denotes a multiindex and $\partial^{\alpha}$ denote partial derivatives with respect to this multiindex. Then the function space $C^m\left(\overline{\Omega}\right)$ with the norm
\begin{equation}
|f|_m:=|f|_{C^m\left(\overline{\Omega}\right) }:=\sum_{|\alpha|\leq m}{\big |}\partial^{\alpha}f{\big |}
\end{equation}
is a Banach space. Here,
\begin{equation}
{\big |}f{\big |}:=\sup_{x\in \Omega}|f(x)|.
\end{equation}

\end{prop}

% Note that a strong spatial scaling $r$ is related to a deviation from a strong semigroup contraction.   However, in refined estimates we mase the idea outlined and local estimates of the form
% \begin{equation}\label{grrh0*}
% \int_{|y|\leq 4\pi \rho r^2 \nu }\frac{lC}{(4 \rho r^2 \nu \tau)^{\delta}|y|^{D-2\delta}}dy\sim(4 \rho r^2 \nu \tau)^{\delta}
% \end{equation} 
% for strong damping parameter $r$ with $4 \rho r^2 \nu\geq 1$. Convolutions with Lipschitz functions then satisfy (let $t_0=0$ for simplicity here)
% \begin{equation}
% {\big |}L\ast G^{\rho,r}_{\nu,i}(\Delta,x){\big |}={\Big |}\int_{}\int_{{\mathbb R}^D}L_x(\sigma,y)G^{\rho,r}_{\nu,i}(\sigma,y)dxd\sigma{\Big |}\leq lC(4\pi \rho r^2 \nu )^{\delta}\Delta^{1-\delta}
% \end{equation}
% for any $\delta \in (0,1)$ and some finite $C>0$.  

Going back to the observation in (\ref{lipob}) note the difference to upper bounds of first order derivatives of the Gaussian, which are integrable only for $\delta >0.5$, due to the Lipschitz continuity of the Leray projection term for strong data in $H^m\cap C^m$ for $m\geq 2$ and the symmetry of the first order spatial derivatives of the Gaussian mentioned above. This upper bound can become smaller than the damping. Indeed for the local solution of the scaled (transformed) equation, i.e., for $v^{\rho,r}_i(\tau,y)=v_i(t,x)$ we get the representation
\begin{equation}\label{Navlerayscheme2irr}
\begin{array}{ll}
 v^{\rho,r,\nu}_i=v^{\rho,r,\nu}_i(t_0,.)\ast_{sp}G^{\rho,r}_{\nu}
-\rho r\sum_{j=1}^D \left( v^{\rho,r,\nu}_jv^{\rho,r,\nu}_i\right) \ast G^{\rho,r}_{\nu,j}\\
\\+\rho r\left( \int_{{\mathbb R}^D}\left( K_D(.-y)\right) \sum_{j,m=1}^D\left( \frac{\partial v^{\rho,r,\nu}_m}{\partial x_j}\frac{\partial v^{\rho,r,\nu}_j}{\partial x_m}\right) (.,y)dy\right) \ast G^{\rho,r}_{\nu,i},
\end{array}
\end{equation}
such that the possible growth  caused by nonlinear terms and estimated by a local time upper bound of the increment of the nonlinear term can be offset by viscosity damping (an argument which holds for positive viscosity).
Let us be a bit more specific.
We search for conditions such that the viscosity damping encoded in the first term on the right side of (\ref{Navlerayscheme2i}) offsets possible growth caused by the nonlinear terms. For this task it is convenient to consider the transformation
\begin{equation}
v^{\nu}_i(t,x)=v^{\rho,r,\nu}_i(\tau,z),~z_i=r x_i,~1\leq i\leq D,~t-t_0=\rho\tau.
\end{equation}
For all $1\leq i,j\leq D$ we have 
\begin{equation}
v^{\nu}_{i,j}(t,x)=v^{\rho,r,\nu}_{i,j}(t,z)r,~v^{\rho,\nu}_{i,j,j}(t,x)=v^{r,\nu}_{i,j,j}(t,z)r^2.
\end{equation}
Hence, under the transformation $z_i=r x_i,~1\leq i\leq D$ the original Cauchy problem for the incompressible Navier Stokes equation (cf. \cite{LL} for the modeling)
\begin{equation}\label{Navlerayequation}
\begin{array}{ll}
 \frac{\partial v^{\nu}_i}{\partial t}-\nu \Delta v^{\nu}_i
\sum_{j=1}^D \left( v^{\nu}_j\frac{\partial v^{\nu}_i}{\partial x_j}\right) =-\nabla_i p,\\
\\
\sum_{i=1}^d\frac{\partial v^{\nu}_i}{\partial x_i}=0,~v^{\nu}_i(0,.)=f_i
\end{array}
\end{equation}
becomes
\begin{equation}\label{Navlerayequationr1}
\begin{array}{ll}
 \frac{\partial v^{\rho,r,\nu}_i}{\partial \tau}-\rho r^2\nu \Delta v^{r,\nu}_i
+\rho r\sum_{j=1}^D \left( v^{r,\nu}_j\frac{\partial v^{r,\nu}_i}{\partial z_j}\right) =-\rho r\nabla_i p^r,\\
\\
r\sum_{i=1}^d\frac{\partial v^{\rho,r,\nu}_i}{\partial z_i}=0,~v^{\rho,r,\nu}_i(0,.)=f_i,
\end{array}
\end{equation}
where for all $t\geq 0$ and $z=rx\in {\mathbb R}^D$ we have
\begin{equation}
p^{\rho,r}(\tau,z)=p(t,x),~p_{,i}(t,x)=p^{\rho,r}_{,i}(\tau,z)r.
\end{equation}
As usual the elimination of the pressure $p^r$ is by application of the divergence operator. From the first equation in (\ref{Navlerayequationr1}) we obtain
\begin{equation}
r^2\sum_{j,i=1}^D  v^{\rho,r,\nu}_{j,i} v^{\rho,r,\nu}_{i,j} =-r^2\Delta  p^{\rho,r},
\end{equation}
which is of the same form as the usual Poisson equation of the non-parametrized velocity, i.e., we have
\begin{equation}
\sum_{j=1}^D  v^{\rho,r,\nu}_{j,i} v^{\rho,r,\nu}_{i,j} =-\Delta  p^{\rho,r}.
\end{equation}
Hence 
\begin{equation}
\int_{{\mathbb R}^D} K_{D,i}(.-y)\sum_{j,m=1}^D\left( \frac{\partial v^{\rho,r,\nu}_m}{\partial x_j}\frac{\partial v^{\rho,r,\nu}_j}{\partial x_m}\right) (.,y)dy
\end{equation}
and the transformed equation becomes
\begin{equation}\label{Navlerayequationr2}
\begin{array}{ll}
 \frac{\partial v^{\rho,r,\nu}_i}{\partial \tau}-\rho r^2\nu \Delta v^{\rho,r,\nu}_i
+\rho r\sum_{j=1}^D \left( v^{\rho,r,\nu}_j\frac{\partial v^{\rho,r,\nu}_i}{\partial x_j}\right) \\
\\-\rho r \int_{{\mathbb R}^D} K_{D,i}(.-y)\sum_{j,m=1}^D\left( \frac{\partial v^{r,\nu}_m}{\partial x_j}\frac{\partial v^{r,\nu}_j}{\partial x_m}\right) (.,y)dy=0.
\end{array}
\end{equation}
Note that the fundamental solution of 
\begin{equation}\label{fundr}
q_{,\tau}-\rho r^2\nu \Delta q=0
\end{equation} 
is explicitly given by 
\begin{equation}
G^{\rho,r}_{\nu}=\frac{1}{\sqrt{4\pi\rho r^2\nu (\tau-\sigma)}^D}\exp\left(-\frac{(x-y)^2}{4\rho r^2\nu (\tau-\sigma)} \right).
\end{equation}
We get
\begin{equation}\label{Navlerayscheme2r}
\begin{array}{ll}
 v^{\rho,r,\nu}_i=v^{\rho,r,\nu}_i(t_0,.)\ast_{sp}G^{\rho,r}_{\nu}
-\rho r\sum_{j=1}^D \left( v^{\rho,r,\nu}_jv^{\rho,r,\nu}_i\right) \ast G^{\rho,r}_{\nu,j}\\
\\+\rho r\left( \int_{{\mathbb R}^D}\left( K_D(.-y)\right) \sum_{j,m=1}^D\left( \frac{\partial v^{\rho,r,\nu}_m}{\partial x_j}\frac{\partial v^{\rho,r,\nu}_j}{\partial x_m}\right) (.,y)dy\right) \ast G^{\rho,r}_{\nu,i}.
\end{array}
\end{equation}
Now the first derivative of the scaled Gaussian $G^{\rho,r}_{\nu}$ is given by
\begin{equation}
G^{\rho,r}_{\nu,i}(\tau,x;\sigma,y)=\left( \frac{-2(x-y)_i}{4\rho r^2\nu (\tau-\sigma)}\right) \frac{1}{\sqrt{4\pi\rho r^2\nu (\tau-\sigma)}^D}\exp\left(-\frac{(x-y)^2}{4\rho r^2\nu (\tau-\sigma)} \right).
\end{equation}
Hence, 
\begin{equation}\label{Navlerayscheme3r}
\begin{array}{ll}
 v^{\rho,r,\nu}_i(\tau,.)=v^{\rho,r,\nu}_i(t_0,.)\ast_{sp}G^{\rho,r}_{\nu}(\tau,.)\\
\\
-\rho r\int_{0}^{\tau}\int_{{\mathbb R}^D}\sum_{j=1}^D \left( v^{\rho,r,\nu}_jv^{\rho,r,\nu}_i\right) (\sigma,y)\left( \frac{-2(.-y)}{4 \nu\rho r^2 (\tau-\sigma)}\right)  G^{\rho,r}_{\nu}\\
\\+\rho r\int_{0}^{\tau}\int_{{\mathbb R}^D}\left( \int_{{\mathbb R}^D}\left( K_D(y-z)\right) \sum_{j,m=1}^D\left( \frac{\partial v^{\rho,r,\nu}_m}{\partial x_j}\frac{\partial v^{\rho,r,\nu}_j}{\partial x_m}\right) (\tau,z)dz\right) \left( \frac{-2(.-y)_i}{4\nu \rho r^2 (\tau-\sigma)}\right)\times\\
\\
\times  G^{\rho,r}_{\nu}(\tau-\sigma,.-y)dyd\sigma.
\end{array}
\end{equation} 
First we consider the damping estimates.  
\begin{itemize}
 \item[i)] First we consider $L^2$-estimates. At time $t_0$ we estimate the growth of the functions  $D^{\beta}_xv_i(t,.),~0\leq |\beta|\leq m,~t\in [t_0,t_0+\Delta]$, where $m\geq 2$ is given, and where
\begin{equation}\label{ass}
\max_{1\leq i\leq n}{\big |}D^{\beta}_xv_i(t_0,.){\big |}_{L^2\cap C}\geq 1.
\end{equation}
If the latter condition is not satisfied for some $\beta$, then there is a $\Delta>0$ such that the respected norm is less equal $1$ for some time $t\in [t_0,t_0+\Delta]$ (because the solutions are locally continuous time curves with values in $H^m\cap C^m$ according to local contraction result), and we need no damping estimate for this part of the $H^m\cap C^m$-norm in the interval $[t_0,t_0+\Delta]$. Proceeding this way we construct an upper bound which is close to the constant $\max_{1\leq i\leq D}{\big |}v_i(t_0,.){\big |}_{H^m\cap C^m}$ times the number of terms in the standard definition of $H^m$-norms. Note that a time interval $t\in [t_0,t_0+\Delta]$ in original time coordinates corresponda to a time interval $\tau\in [t_0,t_0+\Delta_0]$, where for $t=\rho \tau$ we have $\Delta=\rho \Delta_0$.
We apply a spatial Fourier transform with respect to the spatial variables, i.e., the operation
\begin{equation}
{\cal F}(u)(\tau,\xi)=\int_{{\mathbb R}^D}\exp\left(-2\pi i x\xi\right)u(\tau,x)dx,
\end{equation}
in order to analyze the viscosity damping encoded in the first term on the right side of (\ref{Navlerayscheme}) on a time interval $\left[t_0,t_0+\Delta \right]$. For $\tau\in \left[t_0,t_0+\Delta_0 \right]$ and parameters $r,\rho >0$ we have
\begin{equation}
\begin{array}{ll}
{\cal F}\left( v^{\rho ,r,\nu}_i(t_0,.)\ast_{sp}G^{\rho,r}_{\nu}(\tau-t_0,.)\right)={\cal F}\left( v^{\rho,r,\nu}_i(t_0,.)\right) {\cal F}\left( G^{\rho,r}_{\nu}(\tau-t_0,.\right)\\
\\
={\cal F}\left( v^{\rho,r,\nu}_i(t_0,.)\right) \exp\left(-4\pi^2\rho r^2\nu (\tau -t_0) (.)^2 \right),
\end{array}
\end{equation}
where we use that for any $\tau >0$ we have
\begin{equation}
\begin{array}{ll}
{\cal F}\left( G^{\rho,r}_{\nu}(\tau,.)\right) (\tau,\xi)=
{\cal F}\left( \frac{1}{{\sqrt{4\pi \rho r^2 \nu \tau}^D}}
\exp\left(-\frac{(.)^2}{4\nu \rho r^2 \tau} \right)\right)  (\tau,\xi)\\
\\
=\exp\left(-4\pi^2\rho r^2\nu \tau |\xi|^2 \right).
\end{array}
\end{equation}
Let us set $t_0=0$ for simplicity (and without loss of generality of the argument). For $\Delta_0 >0$ small enough (such that, say, $8\pi^2\rho r^2\nu \tau \Delta^2\leq 1$), and for $\tau\in [0,\Delta_0]$ we get
\begin{equation}
\begin{array}{ll}
{\big |}v^{\rho,r,\nu}_i(t_0,.)\ast_{sp}G^{\rho,r}_{\nu}(\tau,.){\big  |}^2_{L^2}\\
\\
=\int_{{\mathbb R}^D}\left( {\cal F}\left( v^{\rho,r,\nu}_i(t_0,.)\right)(\xi) \exp\left(-4\pi^2\rho r^2\nu \tau |\xi|^2 \right)\right) ^2d\xi\\
\\
= \int_{{\mathbb R}^D}\left( {\cal F}\left( v^{\rho,r,\nu}_i(t_0,.)\right)^2(\xi) \exp\left(-8\pi^2\rho r^2\nu \tau |\xi|^2 \right)\right)d\xi\\
\\
=\int_{{\mathbb R}^D\setminus {\{|\xi_j|\leq \Delta_0},1\leq j\leq D\}}\left( {\cal F}\left( v^{\rho,r,\nu}_i(t_0,.)\right)^2(\xi) \exp\left(-8\pi^2\rho r^2\nu \tau |\xi|^2 \right)\right)d\xi\\
\\
+\int_{\{|\xi_j|\leq \Delta_0,~1\leq j\leq D\}}\left( {\cal F}\left( v^{\rho,r,\nu}_i(t_0,.)\right)^2(\xi) \exp\left(-8\pi^2\rho r^2\nu \tau |\xi|^2 \right)\right)d\xi\\
\\
\leq\int_{{\mathbb R}^D}\left( {\cal F}\left( v^{\rho,r,\nu}_i(t_0,.)\right)^2(\xi) \exp\left(-8\pi^2r^2\nu \tau \Delta^2_0 \right)\right)d\xi\\
\\
+{\Big |}\int_{\{|\xi_j|\leq \Delta_0,~1\leq j\leq D\}}{\Big (} {\cal F}\left( v^{\rho,r,\nu}_i(t_0,.)\right)^2(\xi)\times \\
\\
 \times \left( \exp\left(-8\pi^2\rho r^2\nu \tau |\xi|^2 \right)-\exp\left(-8\pi^2\rho r^2\nu \tau \Delta_0^2 \right)\right){\Big )} d\xi{\Big |}\\
\\
\leq {\big |} {\cal F} (v^{\rho,r,\nu}_i)(t_0,.){\big |}_{L^2}^2\exp\left(-8\pi^2\rho r^2\nu \tau \Delta^2_0 \right)
+c^{\Delta}_n
\left(8D\pi^2\rho r^2\nu \tau \Delta_0^{1+D}\right) .
\end{array}
\end{equation}
Here, we use the assumption that $\Delta_0>0$ is small enough (especially $8\pi^2r^2\nu t \Delta_0\leq 1$) and use the abbreviation
\begin{equation}
c^{\Delta}_n:=\sup_{\{|\xi_i|\leq \Delta_0\}}
{\big |}{\cal F}\left( v^{\rho,r,\nu}_i(t_0,.)\right)^2(\xi){\big |}.
\end{equation}
The latter constant is finite (since ${\big |} v^{\nu}_i(t_0,.){\big |}_{H^2\cap C^2}$ is finite, and the square of the latter value is is an upper bound of $c^{\Delta_0}_n$ for sure). 
% Define
% \begin{equation}
% \delta_*v^{\rho,r,\nu}_i(t_0,.)\ast_{sp}G_{\nu}(\tau,.):=v^{\rho,r,\nu}_i(t_0,.)-v^{\rho,r,\nu}_i(t_0,.)\ast_{sp}G_{\nu}(\tau,.)
% \end{equation}
% We get
% \begin{equation}
% \begin{array}{ll}
% {\big |}\delta_*v^{\rho,r,\nu}_i(t_0,.){\big  |}^2_{L^2}={\big |}v^{\rho,r,\nu}_i(t_0,.)\ast_{sp}G^{\rho,r}_{\nu}(\tau,.)-v^{\rho,r,\nu}_i(t_0,.){\big  |}^2_{L^2}\\
% \\
% =\int_{{\mathbb R}^D}\left( {\cal F}\left( v^{\rho,r,\nu}_i(t_0,.)\right)(\xi) \exp\left(-4\pi^2\rho r^2\nu \tau |\xi|^2- {\cal F}\left( v^{\rho,r,\nu}_i(t_0,.)\right)\right)\right) ^2d\xi\\
% \\
% = \int_{{\mathbb R}^D}\left( {\cal F}\left( v^{\rho,r,\nu}_i(t_0,.)\right)^2(\xi) \left( \exp\left(-8\pi^2\rho r^2\nu \tau |\xi|^2 \right)-1\right) \right)d\xi\\
% \\
% \leq {\big |} {\cal F} (v^{\rho,r,\nu}_i)(t_0,.){\big |}_{L^2}^2\exp\left(-8\pi^2\rho r^2\nu \tau \Delta^2 \right)\\
% \\
% +c^{\Delta}_n
% \left(\frac{8}{3}D\pi^2\rho r^2\nu \tau \Delta^{2+D}+\frac{64}{5}D\pi^4\rho r^2\nu^2 \tau^2\Delta^{5+D}\right).
% \end{array}
% \end{equation}

If we take the square root we may use the asymptotics $\sqrt{1+a}=1+\frac{1}{2}a+O(a^2)$.

For $\tau\in [0,\Delta_0]$ and
\begin{equation}
0<\Delta_0 \leq \max\left\lbrace \frac{1}{8\pi^2r^2\nu \max\{c^{\Delta_0}_n,1\}},\frac{1}{2}\right\rbrace 
\end{equation}
we get (the generous) estimate
\begin{equation}\label{vest11}
\begin{array}{ll}
{\big |}v^{\rho,r,\nu}_i(t_0,.)\ast_{sp}G_{\nu}(\tau,.){\big  |}_{L^2}\leq {\big |} {\cal F} (v^{r,\nu}_i)(t_0,.){\big |}_{L^2}\exp\left(-4\pi^2\nu\rho r^2 \tau \Delta_0^2 \right)\\
\\
+c^{\Delta}_n
\left(8D\pi^2\rho r^2\nu \tau \Delta_0^{1+D}\right) \\
\\
\leq {\big |}  v^{r,\nu}_i(t_0,.){\big |}_{L^2}\exp\left(-4\pi^2\nu\rho r^2 \tau \Delta_0^2 \right)+c^{\Delta_0}_n
\left(8D\pi^2\rho r^2\nu \tau \Delta_0^{1+D}\right) .
\end{array}
\end{equation}
If ${\big |}  v^{r,\nu}_i(t_0,.){\big |}_{L^2}$ becomes large or $\Delta >0$ is small enough, then the second summand on right side of (\ref{vest11}) is small compared to the first summand. A similar observation holds for derivatives ${\big |} D^{\beta}_xv^{r,\nu}_i(t_0,.){\big |}_{L^2}$ for $0\leq |\beta |\leq m$.
It is straightforward to obtain analogous estimates for spatial derivatives.
If for $0\leq |\beta|\leq m$ the initial data value ${\big |} {\cal F} (D^{\beta}_xv^{\rho,r,\nu}_i)(t_0,.){\big |}_{L^2}={\big |} D^{\beta}_xv^{\rho,r,\nu}_i)(t_0,.){\big |}_{L^2}$ exceeds a certain level, then for small $\Delta$ the viscosity damping is stronger than possible growth caused by the additional term of order  $\Delta^{D+1}$. As we are interested in the case $D\geq 3$, this is evident, and we may remark in addition that the effect of the additional term becomes smaller as dimension $D$ increases. In items a) and b) below we observe that the damping effect is strong enough in order to offset possible growth caused by the nonlinear terms. Finally note that for $t_0>00$ we get the analogous estimate
\begin{equation}\label{vest111}
\begin{array}{ll}
{\big |}v^{\rho,r,\nu}_i(t_0,.)\ast_{sp}G_{\nu}(\tau,.){\big  |}_{L^2}\leq  {\big |}  v^{r,\nu}_i(t_0,.){\big |}_{L^2}\exp\left(-4\pi^2\nu\rho r^2 (\tau-t_0) \Delta_0^2 \right)\\
\\
+c^{\Delta_0}_n\left(8D\pi^2\rho r^2\nu (\tau -t_0) \Delta_0^{1+D}\right) .
\end{array}
\end{equation}

\item[ii)] Similar considerations hold in other $L^p$-spaces.
\end{itemize}  
Next we show that the viscosity damping can offset possible growth caused by the nonlinear terms. First in item a) we prove the existence of regular upper bounds with respect to ${\big |}.{\big |}_{H^m\cap C^m}$-norms in the case of positive viscosity. Then in item  b) we use  compactness arguments and a similar scheme for a vorticity equation in order to obtain a regular upper bound for a classical solution branch of the Euler equation with respect to the ${\big |}.{\big |}_{C^1\cap H^1}$ norm. Finally, in item c) we obtain upper bounds in stronger norms by auto-controlled schemes, which depend linear on the time horizon.  
Next we introduce Euler Leray data functions of order $l\geq 1$.
\begin{defi}
For $l\geq 1$ and data $\mathbf{g}=(g_1,\cdots ,g_D)^T$ with $g_i\in H^{l+1}\cap C^{l+1}$ for $1\leq i\leq D$ the Euler-Leray data function of order $l$ and of type $1$ is defined by 
\begin{equation}\label{Navlerayfunction}
\begin{array}{ll}
EL^l_1(\mathbf{g}):{\mathbb R}^D\rightarrow {\mathbb R},\\
\\
EL^l_1(\mathbf{g})(y)=\sum_{|\gamma|\leq l}\sum_{j=1}^D D^{\gamma}_x\left( g_jg_{i,j}\right) (y)\\
\\+\sum_{|\gamma|\leq l}\left( \int_{{\mathbb R}^D}\left( K_{D,i}(y-z)\right) \sum_{j,m=1}^DD^{\gamma}_x\left( g_{m,j}g_{j,m}\right) (z)dz\right).
\end{array}
\end{equation}
Furthermore, we define the Euler-Leray data function of type $0$ by
\begin{equation}\label{Navlerayfunction2}
\begin{array}{ll}
EL^l_0(\mathbf{g}):{\mathbb R}^D\rightarrow {\mathbb R},\\
\\
EL^l_0(\mathbf{g})(y)=\sum_{|\gamma|\leq l}\sum_{j=1}^D D^{\gamma}_x\left( g_jg_{i}\right) (y)\\
\\+\sum_{|\gamma|\leq l}\left( \int_{{\mathbb R}^D}\left( K_{D}(y-z)\right) \sum_{j,m=1}^DD^{\gamma}_x\left( g_{m,j}g_{j,m}\right) (z)dz\right),
\end{array}
\end{equation}
\end{defi}
The Euler-Leray functions of order $l\geq 1$ are Lipschitz continuous for regular data $g_i\in H^{l+1}\cap C^{l+1},~1\leq i\leq D$, i.e., for $g_i\in H^{l+1}\cap C^{l+1}$
\begin{equation}
{\big |}EL^l_1(\mathbf{g})(y)-EL^l_1(\mathbf{g})(y'){\big |}\leq L_{l}|y-y'|
\end{equation}
for some finite Lipschitz constant $L_l$. for $g_i\in H^{l+1}\cap C^{l+1}$
\begin{equation}
{\big |}EL^l_0(\mathbf{g})(y)-EL^l_0(\mathbf{g})(y'){\big |}\leq L^0_{l}|y-y'|
\end{equation}
for some finite Lipschitz constant $L_l$. These Lipschitz constants depend on the data norms ${\big |}g_k{\big |}_{H^{l+1}\cap C^{l+1}}$. Nevertheless Lipschitz continuity of the Euler-Leray functions can be computed explicitly for local representations of solutions as in \ref{Navlerayscheme} for strong data due to local contraction results.
 
This Lipschitz continuity can be applied for data at any time $\sigma_0$, i.e., (for fixed argument $x\in {\mathbb R}^D$ and $0\leq |\beta |\leq l$)
\begin{equation}
\begin{array}{ll}
g_i(.)=D^{\beta}_xv^{\rho,r,\nu}_{i} (\sigma_0,x-.)
\end{array}
\end{equation}

\begin{itemize}
\item[a)] 
We consider the case $\nu>0$ and data in $H^m\cap C^m$ at time $t_0\geq 0$. Using the convolution rule from (\ref{Navlerayscheme3r}) we get for all $\tau\in \left[t_0,t_0+\Delta_0 \right]$ and $x\in {\mathbb R}^D$ 
\begin{equation}\label{Navlerayscheme3ar}
\begin{array}{ll}
 v^{\rho,r,\nu}_i(\tau,x)=v^{\rho,r,\nu}_i(t_0,.)\ast_{sp}G^{\rho,r}_{\nu}(\tau-t_0,x)\\
 \\
-\rho r\int_{t_0}^{\tau}\int_{{\mathbb R}^D}\sum_{j=1}^D \left( v^{\rho,r,\nu}_jv^{\rho,r,\nu}_i\right) (\sigma,x-y)\left( \frac{-2(y)_i}{4\nu \rho r^2 (\tau-\sigma)}\right)  G^{\rho,r}_{\nu}(\tau-\sigma,y)dy d\sigma\\
\\+\rho r\int_{t_0}^{\tau}\int_{{\mathbb R}^D}\left( \int_{{\mathbb R}^D}\left( K_D(y-z)\right) \sum_{j,m=1}^D\left( \frac{\partial v^{\rho,r,\nu}_m}{\partial x_j}\frac{\partial v^{\rho,r,\nu}_j}{\partial x_m}\right) (\sigma,x-z)dz\right)\times\\
\\
\times \left( \frac{-2(y)_i}{4 \nu \rho r^2 (\tau-\sigma)}\right)  G^{\rho,r}_{\nu}(\tau-\sigma,y)dyd\sigma.
\end{array}
\end{equation}
For multivariate derivatives of order $m\geq |\beta|=|\gamma|+1,~\gamma_i+1=\beta_i$ we have
\begin{equation}\label{Navlerayscheme4r}
\begin{array}{ll}
 D^{\beta}_xv^{\rho,r,\nu}_i(\tau,x)=D^{\beta}_xv^{\rho,r,\nu}_i(t_0,.)\ast_{sp}G^{\rho,r}_{\nu}\\
 \\
-\rho r\int_{t_0}^{\tau}\int_{{\mathbb R}^D}\sum_{j=1}^D D^{\gamma}_x\left( v^{\rho,r,\nu}_jv^{\rho,r,\nu}_{i,j}\right) (\sigma,x-y)dyd\sigma\\
\\
\left( \frac{-2(y)_i}{4\nu \rho r^2(\tau-\sigma)}\right)  G^{\rho,r}_{\nu}(\tau-\sigma,y)dyd\sigma\\
\\+\rho r\int_{t_0}^{\tau}\int_{{\mathbb R}^D}\left( \int_{{\mathbb R}^D}\left( K_{D,i}(y-z)\right) \sum_{j,m=1}^DD^{\gamma}_x\left( \frac{\partial v^{\rho,r,\nu}_m}{\partial x_j}\frac{\partial v^{\rho,r,\nu}_j}{\partial x_m}\right) (\tau-\sigma,x-z)dz\right) \\
\\
\left( \frac{-2(y)_i}{4 \nu r^2 (\sigma)}\right)  G^{\rho,r}_{\nu }(\tau-\sigma,y)dyd\sigma,
\end{array}
\end{equation}
where in the last step we indicate that we can have the difference $\tau-\sigma$ in the Leray-projection term alternatively. 
Recall from (\ref{grrh0}) that for $\sigma>0$ we have the Gaussian upper bound 
\begin{equation}\label{grrh1}
\forall \delta\in (0,1) ,r,\rho >0~{\big |}G^{\rho,r}_{\nu ,i}(\sigma,y){\big |}\leq \frac{C}{(4 \rho r^2 \nu )^{\delta}\sigma^{\delta}|y|^{D+1-2\delta}}, 
\end{equation}
where (due to Lipschitz continuity of the Leray projection term) $\delta >0$ can be chosen small if this function is convoluted with a Lipschitz continuous function (or a H\"{o}lder continuous function) and $C>0$ is independent of $\rho,r$. There are different estimates for different $\delta \in (0,1)$ as we shall see.
In order to apply Lipschitz continuity of the Euler-Leray function we use local time contraction. 
We have
\begin{lem}\label{contrlem}
Let $t_0\geq 0$, and assume that for some $m\geq 2$ we have
\begin{equation}
{\big |}v^{\rho,r,\nu}_i(t_0,.){\big |}_{H^m\cap C^m}\leq C_0
\end{equation}
For $\delta v^{\rho,r,\nu,k+1}_j=v^{\rho,r,\nu,k+1}_j-v^{r,\nu,k}_j,~1\leq j\leq D$ and $v^{\rho,r,\nu,0}_j=v^{\rho,r,\nu}_j(t_0,.),~1\leq j\leq D$ and $\Delta >0$ small enough we have
\begin{equation}
\sup_{\tau\in [t_0,t_0+\Delta]}{\big |}\delta v^{\rho,r,\nu,k+1}_j(\tau,.){\big |}_{H^m\cap C^m}\leq \frac{1}{2} \sup_{\tau\in [t_0,t_0+\Delta_0]}{\big |}\delta v^{\rho,r,\nu,k}_j(\tau,.){\big |}_{H^m\cap C^m}
\end{equation}
and
\begin{equation}
\sup_{\tau\in [t_0,t_0+\Delta_0]}{\big |}\delta v^{\rho,r,\nu,1}_i(\tau,.){\big |}_{H^m\cap C^m}\leq \frac{1}{2}.
\end{equation}

\end{lem}

For $v^{\rho,r,\nu}_i(\tau,.)\in H^m\cap C^m,~1\leq i\leq D,~\tau\in [t_0,t_0+\Delta_0]$ and all $x\in {\mathbb R}^D$ the function
\begin{equation}\label{lpf}
\begin{array}{ll}
y\rightarrow \int_{{\mathbb R}^D}\left( K_{D,i}(y-z)\right) \sum_{j,k=1}^DD^{\gamma}_x
\left( \frac{\partial v^{\rho,r,\nu}_k}{\partial x_j}\frac{\partial v^{\rho,r,\nu}_j}{\partial x_k}\right) (\tau,x-z)dz,\\
\\
~0\leq |\gamma|\leq m-1
\end{array}
\end{equation}
is Lipschitz continuous with a constant $L_{m-1}>0$ which is independent of $x\in {\mathbb R}^D$. Furthermore as 
\begin{equation}\label{vsup}
\sup_{\tau\in [t_0,t_0+\Delta_0]}{\big |}v^{\rho,r,\nu}_k(\tau,.){\big |}_{H^m\cap C^m}\leq C+1
\end{equation}
the function in (\ref{lpf}) is in $H^{m-1}\cap C^{m-1}$. In the following we shall consider a strong viscosity damping parameter $r$, where
\begin{equation}\label{rhorcond}
4 \rho r^2 \nu\geq 1,
\end{equation}
where we recall that we assumed $\nu >0$ in this item a). 
The function in (\ref{lpf}) is well-defined for all $y\in {\mathbb R}^D$. First we consider the case $|\beta|\geq 1$. From (\ref{Navlerayscheme4r}) (using the convolution rule with respect to time), Lipschitz continuity of the Leray data function with Lipschitz constant $L_m$, and (\ref{grrh1}) we have
 \begin{equation}\label{Navlerayscheme4r*}
\begin{array}{ll}
 {\big |}D^{\beta}_xv^{\rho,r,\nu}_i(\tau,x){\big |}\leq {\big |}D^{\beta}_xv^{\rho,r,\nu}_i(t_0,.)\ast_{sp}G^{\rho,r}_{\nu}(\tau-t_0,x){\big |}\\
 \\
+\rho rL_m\int_{t_0}^{\tau}\int_{B^D}|y|\frac{C}{(4\pi \rho r^2 \nu )^{\delta}(\sigma-t_0)^{\delta}|y|^{D+1-2\delta}}
dyd\sigma +\epsilon,
\end{array}
\end{equation}
where $B^{D}$ is the ball of radius $4 \rho r^2 \nu \geq 1$ (around the origin) , and
\begin{equation}\label{epsilon}
\epsilon =\rho rL_m\int_{t_0}^{t_0+\Delta}\int_{{\mathbb R}^D\setminus B^D}|y|{\big |}G^{r,\rho}_{\nu,i}(\sigma,y){\Big |}dy d\sigma
dyd\sigma 
\end{equation}
becomes small for a small time interval length $\Delta$. Since the upper bound in (\ref{vsup}) holds by the local contraction result, the estimate  in (\ref{Navlerayscheme4r*}) can be transfered to global norms straighforwardly (by application of Young inequalities or by direct estimation convolutions  of the function (\ref{lpf}) with the Gaussian or first order spatial derivatives with the Gaussian). More precisely, 
% since $0\leq |\beta|\leq m$ ${\big |}D^{\beta}_xv^{\rho,r,\nu}_j(\tau,.){\big |}\in C\cap L^2$ we have for $0\leq |\gamma|\leq m-1$ we have ${\big |}v^{\rho,r,\nu}_j(\tau,.){\big |}_{H^{m}\cap C^{m}}\leq C$
% \begin{equation}
% {\Big |}\sum_{j,m=1}^DD^{\gamma}_x\left( \frac{\partial v^{\rho,r,\nu}_m}{\partial x_j}\frac{\partial v^{\rho,r,\nu}_j}{\partial x_m}\right)(\tau,.){\Big |}\leq C_{mD}{\big |}v^{\rho,r,\nu}_m(\tau,.){\big |}_{H^{1}\cap C^1}
% \end{equation}
since for fixed $\tau$ the Gaussian $y\in G^{\rho,r}_{\nu,i }(\tau,y)$ satisfies strong polynomial decay as $|y\uparrow \infty$  (especially has an uppeer bound of form $\frac{c}{1+|y|^{D+2}})$ we have for
\begin{equation}
\begin{array}{ll}
L_x\equiv\int_{{\mathbb R}^D}\left( K_{D,i}(y-z)\right) \sum_{j,k=1}^DD^{\gamma}_x\left( \frac{\partial v^{\rho,r,\nu}_k}{\partial x_j}\frac{\partial v^{\rho,r,\nu}_j}{\partial x_k}\right) (\tau-\sigma,x-z)dz
\end{array}
\end{equation}
for some finite constant $C_{mD}$ and all $0\leq |\gamma|\leq m$ we have
\begin{equation}
\begin{array}{ll}
 {\big |} \rho rD^{\gamma}_xL_{(.)} \ast  G^{\rho,r}_{\nu,i }{\big |}\leq \frac{C_{mD}}{1+|.|^2}\in L^2\cap C,
\end{array}
\end{equation}
where $\ast$ refers to a convolution on the domain $[t_0,t_0+\Delta_0]\times {\mathbb R}^D$.
This implies that the estimates above can be transferred to global norms on small time interval.
Indeed, for $\Delta_0$ small enough and $\tau \in [t_0,t_0+\Delta_0]$ and $4 \rho r^2 \nu\geq 1$ (and $C$ generic) we get
 \begin{equation}\label{Navlerayscheme4r**}
\begin{array}{ll}
 {\big |}D^{\beta}_xv^{\rho,r,\nu}_i(\tau,.){\big |}_{L^2\cap C}\leq {\big |}D^{\beta}_xv^{\rho,r,\nu}_i(t_0,.)\ast_{sp}G^{\rho,r}_{\nu}(\tau,.){\big |}_{L^2\cap C}\\
 \\
+\rho rL_m\int_{t_0}^{\tau}(4\pi \rho r^2 \nu )^{-\delta}C|y|^{2\delta}{\big |}^{4\rho r^2\nu}_0(\sigma-t_0)^{-\delta}d\sigma +\epsilon\\
\\
\leq {\big |}D^{\beta}_xv^{\rho,r,\nu}_i(t_0,.)\ast_{sp}G^{\rho,r}_{\nu}(\tau,.){\big |}_{L^2\cap C}\\
 \\
+\rho rL_m(4 \rho r^2 \nu)^{\delta} \left( \tau-t_0\right) ^{1-\delta} C +\epsilon .
\end{array}
\end{equation}

Recall from (\ref{vest11}) for $t_0=0$ and in general from (\ref{vest111}) for $t_0\geq 0$ and analogous estimates for spatial derivatives that we have a damping estimate
\begin{equation}\label{vest11beta}
\begin{array}{ll}
{\big |}D^{\beta}_xv^{\rho,r,\nu}_i(t_0,.)\ast_{sp}G_{\nu}(\tau,.){\big  |}_{L^2}
\\
\\
\leq {\big |} D^{\beta}_x v^{\rho,r,\nu}_i(t_0,.){\big |}_{L^2}\exp\left(-4\pi^2\nu\rho r^2 \tau \Delta_0^2 \right)+c^{\Delta}_n
\left(8D\pi^2\rho r^2\nu \tau \Delta_0^{1+D}\right),
\end{array}
\end{equation}
which becomes effective for small $\Delta_0 >0$. The last upper bound term in (\ref{Navlerayscheme4r*}) is a constant with respect to space.
Define
\begin{equation}
c^{\Delta}_D:=c^{\Delta}_n
\left(8D\pi^2\rho r^2\nu \tau \Delta_0^{1+D}\right)
\end{equation} 
For $\tau \in [t_0,t_0+\Delta]$ we have 
\begin{equation}\label{Navlerayscheme4r***a}
\begin{array}{ll}
 {\big |}D^{\beta}_xv^{\rho,r,\nu}_i(t_0+\Delta_0,.){\big |}_{L^2}\leq {\big |}  v^{r,\nu}_i(t_0,.){\big |}_{L^2}\exp\left(-4\pi^2\nu\rho r^2 (\tau-t_0) \Delta_0^2 \right)+c^{\Delta}_D\\
 \\
+\rho rL_m(4 \rho r^2 \nu)^{\delta} \left( \Delta_0^{1-\delta}\right) C +\epsilon.
\end{array}
\end{equation}  
We consider the case $t_0=0$ (the estimates for $t_0>0$ are similar). In this case the relation on (\ref{Navlerayscheme4r***a}) shows us that
\begin{equation}
{\big |}D^{\beta}_xv^{\rho,r,\nu}_i(\Delta_0,.){\big |}_{L^2}\leq {\big |}D^{\beta}_xv^{\rho,r,\nu}_i(0,.){\big |}_{L^2}
\end{equation}
if
\begin{equation}\label{ss}
\begin{array}{ll}
{\big |}  v^{\rho, r,\nu}_i(0,.){\big |}_{L^2}\left( \exp\left(-4\pi\nu\rho r^2 \tau \Delta_0^2 \right)-1\right) +c^{\Delta}_D\\
 \\
+\rho rL_m(4 \rho r^2 \nu)^{\delta} \left( \Delta_0^{1-\delta}\right) C +\epsilon \leq 0
\end{array}
\end{equation}
Note that we realize a large damping parameter ($r\geq 1$, w.l.o.g.) in the refined estimate such that
\begin{equation}
4\rho r^2 \nu\geq 1, 
\end{equation}
which imposes the condition
\begin{equation}
r^{1+2\delta}< r^2,~\mbox{ or}~\delta \in \left(0,\frac{1}{2}\right). 
\end{equation}

For a small time interval $\Delta >0$  the positive real number $\epsilon$ in (\ref{epsilon}) becomes arbitrarily small. Furthermore, as $c^{\Delta}_D\downarrow 0$ as  $\Delta_0 \downarrow 0$ with $\Delta_0 ^{D+1}$, and the damping factor  ${\big |}  v^{\rho, r,\nu}_i(0,.){\big |}_{L^2}\left(1- \exp\left(-4\pi^2\nu\rho r^2 \tau \Delta_0^2 \right)\right)$  is dominant for $\tau =\Delta_0$ (recall $t_0=0$) if ${\big |}  v^{\rho, r,\nu}_i(0,.){\big |}_{L^2}$ is greater than a certain threshold (say ${\big |}  v^{\rho, r,\nu}_i(0,.){\big |}_{L^2}=1$) such that $c^{\Delta}_D$ is relatively small compared to the modulus of the main part of this damping term, i.e., small compared to the modulus of \begin{equation}\label{mainpart}
{\big |}  v^{\rho, r,\nu}_i(0,.){\big |}_{L^2}\left( -4\pi^2\nu\rho r^2 \Delta_0 \Delta_0^2 \right).
\end{equation}
Next we consider the conditions such that the modulus of the main damping part is larger than the last term (\ref{Navlerayscheme4r***a}) (the last term for $\epsilon$ which is comparatively small and can be neglected).  Here we observe the exponents of the parameters $\rho,r,\nu,\Delta_0$ in (\ref{mainpart}) compared to the exponents of the parameters $\rho,r,\Delta_0$ of the last term in (\ref{Navlerayscheme4r***a}). For $\tau=\Delta_0$ (case $t_0=0$) in (\ref{mainpart}) we have (for fixed $\nu>0$) the parameter dependence
\begin{equation}\label{depmain}
\sim \nu \rho r^2\Delta_0^3,
\end{equation}
and for the last term in (\ref{Navlerayscheme4r***a}) we have the parameter dependence
\begin{equation}\label{deplast}
(\rho )^{1+\delta}( r )^{1+2\delta}(\nu )^{\delta} \left( \Delta_0^{1-\delta}\right).
\end{equation}
The argument is simplified if we assume $\nu=1$ which is allowed by scaling. However, we shall remark below that a variation of the following argument works for general $\nu >0$.  The viscosity limit is considered in  the next two sections below.
Choosing a small step size parameter $\rho$ and, say $\rho=\Delta_0^{\mu}$ with $\Delta_0<1$ small, we have
\begin{equation}\label{rho}
(\rho )^{1+\delta}( r )^{1+2\delta} \Delta_0^{1-\delta}
= \Delta_0^{\mu(1+\delta)+1-\delta}( r )^{1+2\delta}~~~\mbox{ (case $\nu=1$)},
\end{equation}
and for this choice of parameters the damping term has the dependence (in case $\nu=1$)
\begin{equation}
\rho r^2\Delta_0^3=\Delta_0^{\mu +3}r^2~~\mbox{ (case $\nu=1$)}.
\end{equation}
Hence the damping is stronger then the possible growth of the nonlinear term for ${\big |}D^{\beta}_xv_i(t_0,.){\big |}_{L^2}\geq 1$, for some $0\leq |\beta|\leq m$, $D\geq 3$, and time interval $\Delta_0$ that is small enough, if
\begin{equation}
\mu(1+\delta)+1-\delta>\mu +3~\mbox{and}~\delta \in \left(0,\frac{1}{2}\right), 
\end{equation}
or
\begin{equation}\label{condfin}
\mu >\frac{2+\delta}{\delta}~\mbox{and}~\delta \in \left(0,\frac{1}{2}\right).
\end{equation}
The estimates in the case $|\beta|=0$ are analogous where the Lipschitz constant $L_m$ has to be replaced by $L^0_m$. Note that we can still satisfy \ref{rhorcon}, which becomes
\begin{equation}
4\rho r^2\geq 1\mbox{ for $\nu=1$,}
\end{equation}
as we may choose $r=\frac{1}{\Delta_0^{\frac{\mu}{2}}}$ for example.
\begin{rem}
For arbirary $\nu >0$ we compare the $\nu$-dependent parameter dependence in (\ref{depmain}) and in (\ref{deplast}). Since
\begin{equation}
\nu^{\delta}\leq \nu~~\mbox{if $\nu\geq 1$ and $\delta\in (0,1)$},
\end{equation}
the parameter choice above is suitable also for all finite $\nu\geq 1$.
For $\nu <1$ and still choosing a small step size parameter $\rho$ and, say $\rho=\Delta^{\mu}$, we have to compare the parmeter dependence
\begin{equation}\label{rho}
\nu^{\delta}(\rho )^{1+\delta}( r )^{1+2\delta} \Delta_0^{1-\delta}
= \nu^{\delta}\Delta_0^{\mu(1+\delta)+1-\delta}( r )^{1+2\delta}~~~\mbox{ (case $\nu=1$)},
\end{equation}
with the parameter dependence of the damping term 
\begin{equation}
\nu\rho r^2\Delta_0^3=\Delta_0^{\mu +3}r^2.
\end{equation}
Hence, the damping is stronger then the possible growth of the nonlinear term for ${\big |}D^{\beta}_xv_i(t_0,.){\big |}_{L^2}\geq 1$, for some $0\leq |\beta|\leq m$, $D\geq 3$, and time interval $\Delta_0$ that is small enough, if
\begin{equation}
\mu(1+\delta)+1-\delta>\mu +3~\mbox{and}~\delta \in \left(0,\frac{1}{2}\right), 
\end{equation}
or
\begin{equation}\label{condfin2}
\mu >\frac{2+\delta}{\delta}~\mbox{and}~\delta \in \left(0,\frac{1}{2}\right),
\end{equation}
and
\begin{equation}
\nu^{\delta}r^{1+2\delta}\leq \nu r^2~\mbox{ or }~\nu^{\delta-1}\leq r^{1-2\delta}, i.e.,
\end{equation}
\begin{equation}\label{rchoice}
r\geq \left( \frac{1}{\nu}\right)^{\frac{1-\delta}{1-2\delta}}.
\end{equation}
Note that for the proposed choice (with equality in (\ref{rchoice})) we have
\begin{equation}
4\nu\rho r^2=4\nu \Delta_0^{\mu}\left( \frac{1}{\nu}\right)^{2\frac{1-\delta}{1-2\delta}}=4\nu \Delta_0^{\frac{2+\delta}{\delta}}\left( \frac{1}{\nu}\right)^{2\frac{1-\delta}{1-2\delta}-1}>1
\end{equation}
for $\delta\in (0,0.5)$ close to $0.5$.
 is then realized by the choice. 
\end{rem}

\begin{rem}
For $t_0>0$ analogous estimates hold where $\tau$ is replaced by $\tau -t_0$. It follows that we have constructed global regular upper bounds for multiples of a certain small time interval $\Delta$, where we get a global upper bound for all time from  this by the local time contraction result. This is a global regular upper bound for ${\big |}v^{\rho,r,\nu}_i(t,.){\big |}_{H^m\cap C^m}$ for $m\geq 2$ which transfers to a global regular upper bound for ${\big |}v^{\nu}_i(t,.){\big |}_{H^m\cap C^m}$ for $m\geq 2$ if multiplied by $r^m$ (due to the terms of the highest regularity in the definition of the $H^m\cap C^m$-norm). 
% However, in this case even stronger estimates hold if we choose smaller $\delta$ with $\delta \in \left( 0,\frac{1}{14}\right)$. In this case we get $c^{\Delta}_D\in O(\Delta^{7+\epsilon*)}$ for small $\epsilon*>0$ such that the growth terms integrate to a muber which is compartiavely small comared to the damoing at each step (if the time step size $\Delta >0$ is small enough). Then for $t_0>0$ the damping has a stronger effect. For $\Delta_0=(\tau-t_0)$ we have
% \begin{equation}\label{ss}
% \begin{array}{ll}
% {\big |}  v^{\rho, r,\nu}_i(t_0,.){\big |}_{L^2} \exp\left(-4\pi\nu\rho r^2 (\tau-t_0) \Delta^2_0 \right)+c^{\Delta}_D\\
% \\
% \leq {\big |}  v^{\rho, r,\nu}_i(0,.){\big |}_{L^2}{\Big (}\exp\left(-4\pi\nu\rho r^2 (-t_0) \Delta^2_0 \right)+
% \int_0^{t_0} c^{\Delta}_Dd\tau{\Big )}\times\\
% \\
% \times  \exp\left(-4\pi\nu\rho r^2 (\tau-t_0) \Delta^2_0 \right)
% +c^{\Delta}_D\\
% \\
% \leq  {\big |}  v^{\rho, r,\nu}_i(0,.){\big |}_{L^2}{\Big (} \exp\left(-4\pi\nu\rho r^2 \tau \Delta^2_0 \right)+
% \int_0^{t_0} c^{\Delta}_Dd\tau{\Big )} 
% +c^{\Delta}_D+\epsilon',
% \end{array}
% \end{equation}
% where
% \begin{equation}
%  \epsilon'=\left( \exp\left(-4\pi\nu\rho r^2 (\tau-t_0) \Delta^2_0 \right)-1\right) \left( \int_0^{t_0} c^{\Delta}_Dd\tau{\Big )} \right).
% \end{equation}
% Hence stronger estimates hold for $t_0>0$.
\end{rem}

 \item[b)] The argument of item a) depends essentially on $\nu >0$ as we used 
 \begin{equation}
 4\nu \rho r^2> 1,~\mbox{or }\sqrt{4\nu \rho r^2}> 1.
 \end{equation}
The choice of this 'radius' for the local estimate ensured that the complementary integrals $\int_{{\mathbb R}^D\setminus B_1}\cdots$ become small as $4\nu \rho r^2$ becomes large for fixed $\nu >0$.
Furthermore, since the nonlinear term and the damping term have the same time scaling $\rho$, we have to impose that $r$ depends (somewhat reciprocally) on $\nu$. The main task of this item to get a finite upper bound $C>0$ 
\begin{equation}
{\big |}v^{\nu}_i(\tau,.){\big |}_{H^1\cap C^1}\leq C.
\end{equation}
which is independent of $\nu$ (small). It is clear that the main difficulty for their additional task (compared to item a)) is the estimate for the first order spatial derivatives. One possibility to achieve this task are auto-controlled schemes, which we consider in the next section. However, an alternative approach is to consider the vorticity equation and observe that the it scales similar as the original Navier Stokes equation for the velocity. Indeed, for the same coordinate transformation $(t,x)\rightarrow (\tau,y)$ as above, and for
\begin{equation}
\omega_i^{\rho,r,\nu}(\tau,y):=\omega^{\nu}_i(t,x),
\end{equation}
where from (\ref{vort}) we get the scaled equation
\begin{equation}\label{vort2}
\begin{array}{ll}
 \frac{\partial \omega^{\rho, r,\nu}_i}{\partial \tau}
 -\rho r^2\nu \Delta \omega^{\rho, r,\nu}_i+\rho r\sum_{j=1}^3v^{\rho,r,\nu}_j\frac{\partial \omega^{\rho,r,\nu}_i}{\partial y_j}\\
 \\
=\rho r\sum_{j=1}^3\frac{1}{2}\left(\frac{\partial v^{\rho, r,\nu}_i}{\partial x_j}+\frac{\partial v^{\rho, r,\nu}_j}{\partial x_i} \right)\omega^{\rho, r,\nu}_j.
\end{array}
\end{equation}
Here $v^{\rho,r,\nu}_i,¸1\leq i\leq D$ denotes the scaled original velocity related to a positive viscosity $\nu >0$, where we have
\begin{equation}\label{vel}
v^{\rho,r,\nu}(\tau,y)=\int_{{\mathbb R}^3}K_3(y-z)\omega^{\rho,r,\nu}(\tau,z)dz,~\mbox{where}~K_3(x)h=\frac{1}{4\pi}\frac{x\times h}{|x|^3}.
\end{equation}
We can then argue as in item a) and conclude that we have an upper bound
\begin{equation}
\max_{1\leq i\leq D}\sup_{\tau\in [t_0,t_0+\Delta_0]}{\big |}\omega^{\rho,r,\nu}_i(\tau,.){\big |}_{L^2\cap C}\leq C.
\end{equation}
Then we may use the Bio-Savart law in (\ref{vel}), the Lipschitz continuity and symmetry of the kernel $K_3$  and its first order spatial derivatives $K_{3,j},~1\leq j\leq D$, in order to obtain upper bounds of $v^{\rho,r,\nu}_i(\tau,.),~1\leq i\leq D$ with respect to the $H^1\cap C^1$-norm which are independent of the viscosity $\nu$.  

\item[c)] As an alternative argument to item b) and in order to obtain $\nu$-independent upper bounds in stronger norms we consider auto-controlled schemes. Let us first mention that for auto-controlled schemes which work with small spatial parameter $r>0$ it is sufficient to construct upper bound for the spatial derivatives, i.e., establish 
\begin{equation}\label{beta}
\sup_{\tau \in [0,T]}{\big |}D^{\beta}_xv^{\rho,r}_i(\tau,.){\big |}_{L^2\cap C}\leq C_{\beta}
\end{equation}
 for finite $C_{\beta}$, where $1\leq |\beta|\leq m\geq 2$ using the probabilistic representations in terms of spatial first order derivatives of the Gaussian. Here $\rho$ may depend on the initial data and dimension. The upper bound in the case $\beta=0$ may be obtained the in a seperated argument straightforwardly (using the information in (\ref{beta})). Alternatively, the special structure of the convoluted Leray projection term may be used for $\beta =0$ starting with a scaled version of the representation in (\ref{Navlerayscheme2i}). Derivative rules for convolutions can be used for $v^{\rho,r,\nu}_i(\tau,.)\in H^m\cap C^m$ for $m\geq 2$ such that convolutions of terms of the form $v^{\rho,r,\nu}_{m,j,k}(\tau,.)v^{\rho,r,\nu}_{l}(\tau,.)$ can be estimated similarly as in the case of $|\beta|>0$.
   In item a) and item b) we have considered upper bounds where the condition $\rho r^2\nu \geq 1$ is satisfied. These estimates work in the case of the Navier Stokes equation. However, if we want to extend this argument to the Euler equation, then  this implies that the time step size $\Delta$ of a local iteration scheme has to go to zero as $\nu\downarrow 0$ (i.e., the number of time steps goes to infinity as $\nu$ goes to zero), i.e. we consider analytical limits of the scheme in order to obtain abstract existence results. If we use auto-controlled schemes, then we can eliminate this condition and allow for 
\begin{equation}
\nu\downarrow 0,~\rho r^2\nu\downarrow 0,~\mbox{where }~\rho,r \mbox{ are finite.}
\end{equation} 
Finite $r$ and $\rho$ have the advantage that we can use finite time steps for local iteration steps in the viscosity limit. First we consider the alternative approach of auto-controlled schemes in case of the Navier Stokes equation, i.e., in case of positive viscosity $\nu >0$.

Recall that the first derivative of the scaled Gaussian $G^{\rho,r}_{\nu}$ is given by
\begin{equation}
G^{\rho,r}_{\nu,i}(\tau,x;0,0)=\left( \frac{-2x_i}{4\rho r^2\nu\tau}\right) \frac{1}{\sqrt{4\pi\rho r^2\nu \tau}^D}\exp\left(-\frac{x^2}{4\rho r^2\nu \tau} \right).
\end{equation}
Hence, 
\begin{equation}\label{Navlerayscheme3rc}
\begin{array}{ll}
 v^{\rho,r,\nu}_i(\tau,.)=v^{\rho,r,\nu}_i(t_0,.)\ast_{sp}G^{\rho,r}_{\nu}(\tau,.)\\
\\
-\rho r\int_{0}^{\tau}\int_{{\mathbb R}^D}\sum_{j=1}^D \left( v^{\rho,r,\nu}_jv^{\rho,r,\nu}_i\right) (\sigma,y)\left( \frac{-2(.-y)}{4 \nu\rho r^2 (\tau-\sigma)}\right)  G^{\rho,r}_{\nu}\\
\\+\rho r\int_{0}^{\tau}\int_{{\mathbb R}^D}\left( \int_{{\mathbb R}^D}\left( K_D(y-z)\right) \sum_{j,m=1}^D\left( \frac{\partial v^{\rho,r,\nu}_m}{\partial x_j}\frac{\partial v^{\rho,r,\nu}_j}{\partial x_m}\right) (\tau,z)dz\right)\times \\ 
\\
\times 
\left( \frac{-2(.-y)_i}{4\nu \rho r^2 (\tau-\sigma)}\right) G^{\rho,r}_{\nu}(\tau-\sigma,.-y)dyd\sigma.
\end{array}
\end{equation}

We have constructed the upper bound 
\begin{equation}\label{ordernu}
\sim (\rho )^{1+\delta}( r )^{1+2\delta}(\nu )^{\delta} \left( \Delta_0^{1-\delta}\right).
\end{equation}
for the Leray term
\begin{equation}\label{Navlerayscheme3r*}
\begin{array}{ll}
{\big |}\rho r\int_{t_0}^{t_0+\Delta_0}\int_{{\mathbb R}^D}\left( \int_{{\mathbb R}^D}\left( K_D(y-z)\right) \sum_{j,m=1}^D\left( \frac{\partial v^{\rho,r,\nu}_m}{\partial x_j}\frac{\partial v^{\rho,r,\nu}_j}{\partial x_m}\right) (\tau,z)dz\right)\\
\\
\left( \frac{-2(.-y)_i}{4\nu \rho r^2 (\tau-\sigma)}\right)G^{\rho,r}_{\nu}(\tau-\sigma,.-y)dyds{\big |}_{H^2\cap C^2}.
\end{array}
\end{equation} 
This estimate was constructed in the case $\rho r^2\nu >1$. However as the first derivative of the convoluted Gaussian has the factor $\exp\left(-\frac{x^2}{4\rho r^2\nu (\tau-\sigma)}\right) $ (for $0\leq \tau-\sigma \leq \Delta_0$, upper bounds estimates with parameter dependence as in (\ref{ordernu}) can be obtained for $\rho r^2\nu <1$ for fixed finite $\rho,\nu >0$ if we have the dependence
\begin{equation}
r^2<\Delta_0,~\mbox{or}~\mu <\frac{1}{2}~\mbox{for }~r=\Delta^{\mu},~\mbox{with}~\mu >0.
\end{equation}
The choice of the parameter $\delta$ in the estimate of the Gaussian is also different. recall that convolutions of Lipschitz continuous terms with first order derivatives of the Gaussian can be estimated with $\delta <0.5$ (this does not hold for first order spatial derivatives of the Gaussian itself). 
For  a strong damping parameter $r>1$ we have chosen $\delta \in (0,0.5)$ in order to have a dominant viscosity damping. In auto-controlled schemes we can chose $r>0$ small and $\delta\in (0.5,1)$ alternatively. Note that for fixed $\rho >0$ first order spatial dervatives of the velocity components scale as 
\begin{equation}
v^{\rho,r}_{i,j}(\tau,.)=\frac{v_{i,j}}{r}(t,.),
\end{equation}
hence we have an appropriate scaling of upper bounds for the Leray projection increments of order
\begin{equation}
\sim r^{1+2\delta}(\nu )^{\delta} \left( \Delta_0^{1-\delta}\right)\lesssim \Delta_0^{1+\epsilon}
\end{equation}
for some $\epsilon >0$ if 
\begin{equation}
\mu(1+2\delta )>\delta~\mbox{or}~\mu>\frac{\delta}{1+2\delta},
\end{equation}
and this is perfectly consistent with the condition
\begin{equation}
\mu<\frac{1}{2}
\end{equation}
above which we needed in order to have in order to transfer our upper bound estimates abive to the case where $\rho r^2\nu <1$. 
 Note that the increment of the potential damping term (cf. of multivariate spatial derivatives $D^{\beta}_xv^{\rho,r}_i,~0\leq |\beta|\leq 2$ is of the form
\begin{equation}
\sim \Delta_0 D^{\beta}_xv^{\rho,r}_i=\Delta_0 \frac{D^{\beta}_xv_i}{r}.
\end{equation}
which offsets possible growth cause by the Leray projection term increment for $\Delta_0\in (0,1)$ small enough and if 
$D^{\beta}_xv^{\rho,r}_i(t_0,.)$ exceeds a certain threshold. This argument holds for positive viscosity $\nu>0$, i.e., for the Navier Stokes equation, but it is not suitable for the viscosity limit $\nu \downarrow 0$, i.e., for the Euler equation.

For the Euler equation the situation is less comfortable.
As $\rho r^2\nu\downarrow 0$ this alternative agument does not depend on the viscosity damping. Hence we may assume that $r>0$ is small. Furthermore, we use a time scaling works for $\rho>0$, i.e., where $\rho$ depends only on dimansion and the initial data. 
Given $t_0\geq 0$ and a time interval $\left[t_0,t_0+\Delta_0 \right]$ for some $\Delta_0 \in (0,1)$ we consider a comparison function $u^{\rho,r,\nu,t_0}_i,~1\leq i\leq D$ by a global time transformation with local time dilatation on a short time interval, i.e., we consider the transformation
\begin{equation}\label{utr1}
\begin{array}{ll}
(1+\tau )u^{\rho,r,\nu,t_0}_i(s,y)=v^{\rho,r,\nu}_i(t,y),\\
\\
\mbox{where}~s=\frac{\tau-t_0}{\sqrt{1-(\tau-t_0)^2}},~\tau-t_0\in [0,\Delta_0],~1\leq i\leq D
\end{array}
\end{equation}
where $\Delta_0 \in [0,1)$. 
Using the abbreviation $\Delta \tau=\tau-t_0$, we have
\begin{equation}\label{utr2}
 v^{\rho,r,\nu}_{i,t}=
u^{\rho,r,\nu,t_0}_i
+(1+\tau )u^{\rho,r,\nu,t_0}_{i,s}\frac{ds}{dt},~\frac{ds}{dt}
=\frac{1}{\sqrt{1-\Delta \tau^2}^3}.
\end{equation}
We choose $[0,\Delta]=[0,0.5]$ for convenience. The equation for $u^{\rho,r,\nu,t_0}_i,~1\leq i\leq D$ becomes
\begin{equation}
\begin{array}{ll}\label{utr3}
u^{\rho,r,\nu,t_0}_{i,s}+\frac{\sqrt{1-\Delta \tau^2}^3}{1+\tau} r^2\nu \Delta u^{\rho,r,\nu,t_0}_i
+ r\frac{\sqrt{1-\Delta \tau^2}^3}{1+\tau}\sum_{j=1}^D \left( u^{\rho,r,\nu,t_0}_j\frac{\partial u^{\rho,r,\nu,t_0}_i}{\partial x_j}\right) \\
\\
-r\frac{\sqrt{1-\Delta \tau^2}^3}{1+\tau} \int_{{\mathbb R}^D} K_{D,i}(.-y)\sum_{j,m=1}^D\left( \frac{\partial u^{\rho,r,\nu,t_0}_m}{\partial x_j}\frac{\partial u^{\rho,r,\nu,t_0}_j}{\partial x_m}\right) (.,y)dy\\
\\
-\frac{\sqrt{1-\Delta \tau^2}^3}{1+\tau}u^{r,\nu,t_0}_i=0,\\
\\
~u^{\rho,r,\nu}(0,.)=v^{\rho, r,\nu}(t_0,.).
\end{array}
\end{equation}
Now assume that an arbitrary large but fixed time horizon $T>0$ is given.
Assume that for $m\geq 2$ and for some $0\leq t_0<T$ we have a global upper bound
\begin{equation}
\max_{1\leq i\leq D}{\Big |}v^{\rho,r,\nu}_i(t_0,.){\Big |}_{H^m\cap C^m}\leq C(1+t_0)<\infty,
\end{equation}
where $C>0$ is independent of $\nu$. We may assume that $C\geq 3$. Then 
\begin{equation}
\max_{1\leq i\leq D}{\Big |}u^{\rho,r,\nu,t_0}_i(0,.){\Big |}_{H^m\cap C^m}=\frac{1}{1+t_0} \max_{1\leq i\leq D}{\Big |}v^{\rho,r,\nu}_i(t_0,.){\Big |}_{H^m\cap C^m}\leq C.
\end{equation}
 Using local iteration schemes and representations as in (\ref{Navlerayscheme3ar}) and (\ref{Navlerayscheme4r}), local upper bounds for Gaussians and for first order spatial derivatives of the Gaussian, Young inequalities, Fourier transforms and Lipschitz continuity of the Euler-Leray function, it is straightforward to obtain a local time contraction result in Lemma \ref{contrlem*} below.
 We conclude that for $r$ as in (\ref{contrlem*}) and for $C\geq 3$ we have
\begin{equation}
\max_{1\leq i\leq D}{\Big |}u^{\rho,r,\nu,t_0}_i(\Delta_0,.){\Big |}_{H^m\cap C^m}\leq C.
\end{equation} 
This implies that 
\begin{equation}
\begin{array}{ll}
\max_{1\leq i\leq D}{\Big |}v^{r,\nu}_i(t_0+\Delta_0,.){\Big |}_{H^m\cap C^m}\\ 
\\
\leq \max_{1\leq i\leq D}(1+t_0+\Delta_0){\Big |}u^{r,\nu,t_0}_i(\Delta_0,.){\Big |}_{H^m\cap C^m}
 \leq C(1+t_0+\Delta_0).
\end{array}
\end{equation} 
Here $r$ depends only some time invariant constants (cf. (\ref{contrlem*}) below) such that for all integers $k\geq 0$ with $k\Delta \leq T$ we have 
\begin{equation}
\begin{array}{ll}
\max_{1\leq i\leq D}{\Big |}v^{\rho,r,\nu}_i(k\Delta_0,.){\Big |}_{H^m\cap C^m}\leq \max_{1\leq i\leq D}(1+k\Delta_0){\Big |}u^{\rho,r,\nu,t_0}_i(0,.){\Big |}_{H^m\cap C^m}\\
\\
\leq C(1+k\Delta_0).
\end{array}
\end{equation}
Then by a local contraction result similar as in (\ref{contrlem}) above we can interpolate and have
\begin{equation}
\begin{array}{ll}
\max_{1\leq i\leq D}{\Big |}v^{\rho,r,\nu}_i(\tau,.){\Big |}_{H^m\cap C^m}\leq \max_{1\leq i\leq D}(1+\tau){\Big |}u^{\rho,r,\nu,t_0}_i(0,.){\Big |}_{H^m\cap C^m} \\
\\
\leq 2C(1+\tau).
\end{array}
\end{equation}
The local time contraction result used in this argument is
\begin{lem}\label{contrlem*}
Let an arbitrarily large but finite time horizon $T>0$ be given and let $C_0\geq 3$ and $\Delta_0 =0.5$. Let $t_0\geq 0$ and assume that for some $m\geq 2$ we have
\begin{equation}
\max_{1\leq i\leq D}{\big |}u^{\rho,r,\nu,t_0}_i(0,.){\big |}_{H^m\cap C^m}\leq C_0.
\end{equation}
Recall that $D$ denotes the dimension and define the constants 
\begin{equation}
C_K=\max\left\lbrace {1,\big |}K_{D,i}{\big |}_{L^1(B_1)}+{\big |}K_{D,i}{\big |}_{L^2({\mathbb R}^D\setminus B_1)}\right\rbrace ,
\end{equation}
where $B_1$ is the ball of radius $1$ in ${\mathbb R}^D$, and
\begin{equation}
\begin{array}{ll}
C_G=1+ {\big |}G^{\rho,r}_{\nu}{\big |}_{L^1((0,T)\times B_1)}+{\big |}G^r_{\nu}{\big |}_{L^2\left((0,T)\times {\mathbb R}^D\setminus B_1\right) }\\
\\
+\max_{1\leq i\leq D}{\big |}G^r_{\nu,i}{\big |}_{L^1((0,T)\times B_1)}+{\big |}G^{\rho,r}_{\nu,i}{\big |}_{L^2\left((0,T)\times {\mathbb R}^D\setminus B_1\right) }.
\end{array}
\end{equation}
For $\delta u^{\rho r,\nu,t_0,k+1}_j=u^{\rho,r,\nu,t_0,k+1}_j-u^{\rho,r,\nu,t_0k}_j,~1\leq j\leq D$ and $u^{\rho,r,\nu,t_0,0}_j=u^{\rho,r,\nu,t_0}_j(0,.),~1\leq j\leq D$ and $\Delta_0\leq 0.5$ we have for
\begin{equation}
r=\frac{1}{2C_0},~\rho=\frac{1}{2DC_0(C_G+DC_KC_G)}
\end{equation}
\begin{equation}
\sup_{s\in [t_0,t_0+\Delta_0]}{\big |}\delta u^{r,\nu,t_0,k+1}_j(s,.){\big |}_{H^m\cap C^m}\leq \frac{1}{2} \sup_{\tau\in [t_0,t_0+\Delta_0]}{\big |}\delta u^{r,\nu,t_0,k}_j(s,.){\big |}_{H^m\cap C^m}
\end{equation}
and
\begin{equation}
\sup_{s\in [t_0,t_0+\Delta]}{\big |}\delta u^{r,\nu,1}_i(s,.){\big |}_{H^m\cap C^m}\leq \frac{1}{2}.
\end{equation}
Moreover,
\begin{equation}
\max_{1\leq i\leq D}{\big |}u^{r,\nu,t_0}_i(\Delta_0,.){\big |}_{H^m\cap C^m}\leq \max_{1\leq i\leq D}{\big |}u^{r,\nu,t_0}_i(0,.){\big |}_{H^m\cap C^m}\leq C_0.
\end{equation}
Here, $C_0$ can be chosen independently of $\nu >0$.
\end{lem}
 
\end{itemize}

\section{Short-and Long time singularities of the Navier Stokes equations with time dependent force term}
We consider the argument for a global upper bound $v^{\rho,r,\nu}_i(\tau,.)$ of item b) in more detail. Note that our construction of global upper bounds for $v^{\rho,r,\nu}_i,¸1\leq i\leq D$ enforces 
\begin{equation}
r\uparrow \infty \mbox{ as }\nu \downarrow 0~\mbox{ (for uniform $H^1 \cap C^1$ upper bounds)}.
\end{equation}
This is a difference to autocontrlled schemes of item c) where
\begin{equation}
r\downarrow 0 \mbox{ as }\nu \downarrow 0~\mbox{ (for time-linear $H^m \cap C^m$ upper bounds for $m\geq 2$)}.
\end{equation}
Hence a direct upper bound argument for the velocity components as in item b) is not possible for ${\big |}v^{\rho,r,\nu}_i(\tau,.){\big |}_{H^1\cap C^1}$, and we shall observe in the next section that it is possible for  ${\big |}v^{\rho,r,\nu}_i(\tau,.){\big |}_{L^2\cap C^0}$ if we have a strong spatial decay of the data such that strong compactness arguments are available. However, in the special case of the Euler equation we can improve the siuation by consideration of the vorticity equation.
An analogous $L^2$-argument as in (\ref{Navlerayscheme4r***2}) for the vorticity shows that for $\Delta_0 >0$ small enough and data in $H^{m}\cap C^{m},~m\geq 2$ possible growth caused by the nonlinear terms is offset by viscosity damping in the sense that (e.g.) for $\delta \in \left(0, \frac{1}{2}\right) $, a small time interval $[t_0,t_0+\Delta_0]$ , small $\rho\sim \Delta^{\mu} $ with $(\mu,\delta)$ as (\ref{condfin})  we get for all $0\leq |\beta|\leq m$ and a finite constant $C_{\omega}$ the analogous vorticity estimate
 \begin{equation}\label{Navlerayscheme4r***2}
\begin{array}{ll}
 {\big |}D^{\beta}_x\omega^{\rho,r,\nu}_i(t_0+\Delta_0,.){\big |}_{L^2\cap C}\leq {\big |}  \omega^{r,\nu}_i(t_0,.){\big |}_{L^2}\exp\left(-4\pi\nu\rho r^2 (\tau -t_0)\Delta^2 \right)+c^{\Delta}_n\\
 \\
+\rho rL_m(4 \rho r^2 \nu)^{\delta} \left( \Delta_0^{1-\delta}\right) C_{\omega} +\epsilon
\leq {\big |}D^{\beta}_x\omega^{\rho,r,\nu}_i(t_0,.){\big |}_{L^2\cap C},
\end{array}
\end{equation}  
if the latter term is larger than a certain threshold, say, it satisfies
\begin{equation}
{\big |}D^{\beta}_x\omega^{\rho,r,\nu}_i(t_0,.){\big |}_{L^2\cap C}\geq 1.
\end{equation}
Recall that (\ref{Navlerayscheme4r***2}) holds if
\begin{equation}
4 \rho r^2 \nu >\sqrt{4 \rho r^2 \nu}>1.
\end{equation}

We have not observed explicitly that these upper bounds are independent of $\nu$ essentially, in the sense that for $\nu$ small, e.g., $\nu\sim \sqrt{\Delta_0}$ as the time interval becomes small, we can set up a scheme of step size $\Delta_0^{\mu}$ solving for $\omega^{\rho,r,\nu}_i,~1\leq i\leq D$ which has a global regular upper bound which is independent of the size of $\Delta_0$, and, therefore, independent of $\nu$.   
Let us choose $\nu=\sqrt{\Delta_0}$ for small $\Delta_0 >0$, and 
\begin{equation}\label{rchoice}
1\geq r=\frac{1}{\sqrt{\Delta_0}^{\mu+1}},~\rho=\Delta_0^{\mu},~\mu >\frac{2+\delta}{\delta}~\mbox{and}~\delta \in \left(0,\frac{1}{2}\right).
\end{equation}
Then we have
\begin{equation}
4\rho\nu r^2=4 >1,
\end{equation}
where the damping estimate has a stronger scaling with respect to $r$ than the upper bound of the nonlinear term in the sense that 
\begin{equation}
r^2>r^{1+2\delta}~\mbox{for $r> 1$.}
\end{equation}
The growth caused by the nonlinear terms in an interval $[t_0,t_0+\Delta_0]$ is offset by the damping for small $\Delta_0$ under the condition (\ref{rchoice}) if 
\begin{equation}
 \rho r(4 \rho r^2 \nu)^{\delta} \left( \Delta_0^{1-\delta}\right) \leq \nu\rho r^2\Delta_0^3,
\end{equation}
or
\begin{equation}
 \Delta_0^{\mu} \frac{1}{\sqrt{\Delta_0}^{\mu+1}}\left( 4 \Delta_0^{\mu} \frac{1}{\Delta_0^{\mu+1}}\sqrt{\Delta_0}\right) ^{\delta} \left( \Delta_0^{1-\delta}\right) \leq \sqrt{\Delta_0}\Delta_0^{\mu} \frac{1}{\Delta_0^{\mu+1}}\Delta_0^3,
\end{equation}
which means (for small $\Delta_0$)
\begin{equation}
\Delta_0^{\mu-\frac{\mu+1}{2}}\left( 4 \sqrt{\Delta_0}\right) ^{-\delta} \left( \Delta_0^{1-\delta}\right)\leq \Delta_0^{2.5}~\mbox{,or}~\mu \geq 5+3\delta . 
\end{equation}
Hence we the upper bound estimates in item b) are essentially independent of $\nu$ for ${\big |}\omega^{\rho,r,\nu}_i(\tau,.){\big |}_{L^2\cap C}$ and hence for ${\big |}v^{\rho,r,\nu}_i(\tau,.){\big |}_{H^1\cap C^1}$
\begin{rem} We note that Lipschitz continuity or strong H\"{o}lder continuity of the Leray data functions can used in order to obtain local time contraction results based on convolution with first order derivatives of the Gaussian in the viscosity limit. 
More precisely, note that in general iterative solution of an equation
\begin{equation}\label{exeq}
v_{,t}-\nu \Delta v+F(v,\nabla v)=0,~v(t_0,.)=v_{t_0}
\end{equation}
by an iterative convolution scheme of the form 
\begin{equation}\label{fconv}
v^{k+1}=v^k\ast_{sp}\ast G_{\nu}+F(v^k)\ast G_{\nu,j}
\end{equation}
causes principle problems in the limit $\nu\downarrow 0$ (even if considered on a compact domain $\Omega$) due to
a degeneracy (with respect to $\nu$) of a typical upper bounds (for some $\delta \in (0,1)$) obtained, e.g., from intergals of the upper bound
\begin{equation}\label{grrh12}
l|y|{\big |}G^{\rho,r}_{\nu ,i}(\sigma,y){\big |}\leq l\frac{C}{(4\pi \rho r^2 \nu )^{\delta}\sigma^{\delta}|y|^{D-2\delta}}, 
\end{equation}
on a ball of radius $4\pi\rho r^2\nu\geq 1$. However if the analysis of (\ref{fconv}) leads to a functional series $(v^{\nu_k})_{\nu_k}$ for a sequence $(\nu_k)\downarrow 0$ as $k\uparrow \infty$ such that $v^{\nu_k}$ solves (\ref{exeq})  for $\nu_k$ on a time interval $[t_0,t_0+\Delta]$ and we have an upper bound
\begin{equation}
\sup_{t\in [t_0,t_0+\Delta]}{\big |}v^{\nu_k}(t,.){\big |}_{H^m_0(\Omega)}\leq C
\end{equation}
for a strong $H^m_0$-norm (where  $H^m_0\left(\Omega\right)$ is the closure of $C^{\infty}_c(\Omega)$ in $H^m$, and for a finite constant $C>0$ which is independent of $\nu_k$ for all $k\geq 1$), then Rellich's compactness result gives a subsequence $(v^{\nu'_k})_{\nu'_k>0}$ with $\nu'_k\downarrow 0$ such that the limit is in the slightly weaker space $v^{0}(t,.)\in H^{m-\epsilon}$ for any small $\epsilon>0$. For large $m$ (dependent on dimension, of course) the limit function $v^0$ solves the equation in (\ref{exeq}) pointwise for strong data in $H^m_0(\Omega)$. Such a reasoning can be transferred to infinite domains if we have strong polynomial decay of the data and the operator function $F$ preserves this strong polynomial decay. The latter property holds for  the Burgers and Leray projection term for strong polynomial decay as we observe in the next sections.   
\end{rem}
The existence of singular solutions of the incompressible Euler equation implies the existence of singular solutions of the incompressible Navier Stokes equation with time dependent force terms. A classical solution 
\begin{equation}
v_i,~1\leq i\leq D,~v_i(t,.)\in H^m\cap C^m,~t\in [0,T),~m\geq 2
\end{equation}
of the incompressible Euler equation (on the interval $[0,T)$) with a blow-up of vorticity at time $T$ satisfies the incompressible Navier Stokes equation
 \begin{equation}\label{vorticitynav}
\frac{\partial v_i}{\partial t}-\nu \Delta v_i+
\sum_{j=1}^nv_jv_{i,j}-\int_{{\mathbb R}^D} 
K_{D,i}(.-z) \sum_{j,m=1}^D\left( \frac{\partial v^{\nu}_m}{\partial x_j}\frac{\partial v^{\nu}_j}{\partial x_m}\right) (.,z)dz=F_i, 
\end{equation}
with force term $F=\left(F_1,F_2,\cdots ,F_D \right)^T$ (on the same time interval) if
\begin{equation}\label{Fterm}
F_i=-\nu \Delta v_i\in H^{m-2}\cap C^{m-2},~1\leq i\leq D.
\end{equation}
The analysis of global regular solution branches and time reversed Euler equations below shows that $F_i$ is also $L^2$ with respect to time on the time interval $[0,T]$, where a possible large time $T>0$ is the time where the vorticity of the Euler equation blows up. As a consequence of this analysis we note that there are long time singularities.  

\begin{thm}
The Navier Stokes equation  with initial data $h_i\in H^{m}\cap C^{m},~m\geq 2,~1\leq i\leq D$ of strong polynomial decay, i.e., $h_i\in {\cal C}^{m(D+1)}_{pol,m},~1\leq i\leq D$  (cf. the definition of the latter function space in the next section below), and external time-dependent forces $F_i\in H^{m-2}\cap C^{m-2},~1\leq i\leq D$ is not well -posed in general. 
More precisely, for the Navier Stokes equation Cauchy problem in (\ref{vorticitynav}) and for any time $T>0$ there exist time dependent force terms
\begin{equation}
F_i(t,.)\in H^m\cap C^m,~k\geq 0,~\tau\in[0,T]
\end{equation}
and data  $v_i(0,.)\in H^{m}\cap C^{m}\cap {\cal C}^{m(D+1)}_{pol,m} ,~1\leq i\leq 3,~m\geq 2$ such that a regular classical solution of the Navier Stokes equation on the time interval $[0,T)$ has a blow-up or a kink of order $k\geq 1$ at time $T$.
\end{thm}

\section{Global solution branches of the Euler equation}
We have found upper bounds which are independent of the viscosity constant $\nu >0$. These upper bounds are obtained either in $H^1\cap C^1$-space by Lipschitz continuity of the vorticity (strategy of item b) above, or, for $H^m\cap C^m$-spaces, by auto-controlled schemes (strategy of item c) above. In the latter case (strategy of item c)) the upper bounds for spatial derivatives of the velocity component functions depend linearly on the time horizon $T>0$ (so the result of item b) is slightly stronger with respect to upper bounds in ${\big |}.{\big |}_{H^1\cap C^1}$-norms). Strong spatial polynomial decay implies that the Navier-Stokes and Euler equation Cauchy problem can be transformed to a problem with spatially compact domain (with natural Dirichlet boundary conditions at the spatial boundaries). This has the advantage that stronger Banach spaces -$C^1$-Banach spaces as in Proposition 1.1.- can be used in order to construct viscosity limits (instead of Rellich embedding theorems where a small loss of regularity is implied). Next we consider these ideas in more detail.       
In any case it is implied that for any finite $T>0$ there exists a finite constant $\tilde{C}>0$ (which may depend linearly on the time horizon $T$) such that  we have
\begin{equation}\label{upperbound2}
\begin{array}{ll}
\max_{1\leq i\leq D}\inf_{\nu>0}\sup_{t\in [0,T]}{\big |}v^{\nu}_i(t,.){\big |}_{H^m\cap C^m}\leq \tilde{C}.
\end{array}
\end{equation}
As we have an unbounded domain we have to be a little careful concerning compactness arguments. However, the schemes used here preserve some degree of strong polynomial spatial decay. This is not well-known as the convolution effect with the Gaussian decreases the degree of polynomial decay. Therefore we repeat the main arguments for that observation here (cf. remark \ref{rempdl}). If for fixed $\rho,r>0$ a solution $v^{\rho,r,\nu}_i,~1\leq i\leq D$ has polynomial spatial decay of order $m(D+1)-1$ (we are concerned with $D\geq 3$), then standard compactness arguments can be applied.
First we observe the inheritance of polynomial decay of order $m(D+1)-1$ by the scheme if data have spatial polynomial decay of order $m(D+1)$. We define a related function space.
\begin{defi}
For $l\geq 1$ and $m\geq 2$ we define a space of functions which satisfy polynomial decay of order $l\geq 1$ at spatial infinity for multivariate spatial derivatives up to order $m$. More precisely, we define
\begin{equation}
{\cal C}^{l}_{pol,m}={\Big \{} f:{\mathbb R}^D\rightarrow {\mathbb R}: \\
\\
\exists c>0~\forall |x|\geq 1~\forall 0\leq |\gamma|\leq m~{\big |}D^{\gamma}_xf(x){\big |}\leq \frac{c}{1+|x|^l} {\Big \}}. 
\end{equation}
\end{defi}
First we recall -given data of strong spatial polynomial decay- the inheritance of some order of polynomial decay by the local scheme (next remark). If we avoid the convolution of the initial data in the local scheme, then the corresponding local scheme preserves even a stronger form of polynomial decay. First we observe
\begin{rem}\label{rempdl}
Consider first local time  iterative solution schemes $v^{\rho,r,\nu,(k)}_i,~1\leq i\leq n,~k\geq 0$ of the scaled Navier Stokes equation in terms of classical representations using convolutions with the fundamental solution $G^{\rho,r}_{\nu}$  of the equation $\frac{\partial u}{\partial t}-\rho r^2\nu\sum_{j=1}^n \frac{\partial^2 u}{\partial x_j^2}=0$. Let $t_0\geq 0$ and consider a time interval $[t_0,t_0+\Delta]$ for some small $\Delta >0$. Assume that the data at time $t_0$ are of strong polynomial spatial decay, i.e., they satisfy
\begin{equation}
v^{\rho,r,\nu}_i(t_0,.)\in {\cal C}^{m(D+1)}_{pol,m}~\mbox{ for some } m\geq 2.
\end{equation}
Define for all $\tau\geq t_0$ and $x\in {\mathbb R}^n$
\begin{equation}
v^{\rho,r,\nu (0)}_i(\tau,x):=\int_{{\mathbb R}^n}v^{\rho,r,\nu}_i(t_0,.)_i(y) G_{\rho,r,\nu}(\tau,x;0,y)dy=v^{\rho,r,\nu}_i(t_0,.)_i\ast_{sp}G_{\nu},
\end{equation}
and for $k\geq 1$
\begin{equation}\label{Navleraysol}
\begin{array}{ll}
v^{\rho,r,\nu,(k)}_i
=v^{\rho,r,\nu}_i(t_0,.)_i\ast_{sp}G_{\nu}+\sum_{j=1}^n \left( v^{\rho,r,\nu,(k-1)}_j\frac{\partial v^{\rho,r,\nu,(k-1)}_i}{\partial x_j}\right) \ast G_{\nu}\\
\\+\sum_{j,m=1}^n\int_{{\mathbb R}^n}\left( \frac{\partial}{\partial x_i}K_n(.-y)\right) \sum_{j,m=1}^n\left( \frac{\partial v^{\rho,r,\nu,(k-1)}_m}{\partial x_j}\frac{\partial v^{\rho,r,\nu,(k-1)}_j}{\partial x_m}\right) (.,y)dy\ast G_{\nu},
\end{array}
\end{equation}
where we recall that $\ast$ denotes convolution with respect to space and time. In the computation of the increment
\begin{equation}
\delta v^{\rho,r,\nu}_i:=v^{\rho,r,\nu}_i-v^{\rho,r,\nu}_i(t_0,.)\ast_{sp}G_{\nu}
\end{equation}
the convolutive effect of the Gaussian can be eliminated by the consideration of a related iterative solution scheme. Define
\begin{equation}
\delta v^{\rho,r,\nu,(0)}_i(t,x):=0,
\end{equation}
and for $k\geq 1$ define
\begin{equation}\label{Navleraysol}
\begin{array}{ll}
\delta v^{\rho,r,\nu,(k)}_i:=v^{\rho,r,\nu,(k)}_i-v^{\rho,r,\nu,(k-1)}_i.
\end{array}
\end{equation}
At each iteration step $k$ in the contribution of the nonlinear quadratic terms the decreasing order of decay at spatial infinity of each factor caused by convolutive effects is offset by the multiplicative effect in the nonlinear term. It is straightforward to prove that for local time $t\geq 0$ and $m\geq 2$
\begin{equation}
\delta v^{\rho,r,\nu,(k)}(t,.)\in {\cal C}^{m(n+1)}_{pol,m},~v^{\rho,r,\nu,(k)}(t,.)\in {\cal C}^{m(n+1)-\mu}_{pol,m}
\end{equation}
for some $\mu\in (0,1)$ which accounts for the convolution effects of the linear term (in  this case the convoluted initial data).  This leads to a local time solution representation 
\begin{equation}
v^{\rho,r,\nu}_i:=v^{\rho,r,\nu,(0)}_i+\sum_{k\geq 1}\delta v^{\rho,r,\nu,(k)}_i=v^{\rho,r,\nu}_i(t_0,.)\ast_{sp}G_{\nu}+\sum_{k\geq 1}\delta v^{(k)}_i\in {\cal C}^{m(n+1)-\mu}_{pol,m},
\end{equation}
where we have a small loss of spatial decay of the solution compared to the initial data for the local scheme.
Next we observe that this small loss of spatial decay is not increased by the global scheme. Assume that
\begin{equation}
v^{\rho,r,\nu}_i(t_1,.)=\left( v^{\rho,r,\nu}_i(t_0,.)\ast_{sp}G^{\rho,r}_{\nu}\right)(t_1,.)+\left( \sum_{k\geq 1}\delta v^{\rho,r,\nu,(k)}_i\right)(t_1,.)\in {\cal C}^{m(n+1)-\mu}_{pol,m},
\end{equation} 
has been proved for some $t_1>0$. Furthermore assume that
\begin{equation}\label{assdelta}
\sum_{k\geq 1}\delta v^{\rho,r,\nu,(k)}_i(t_0,.)\in {\cal C}^{m(n+1)}_{pol,m}
\end{equation}
We then define a local time iteration scheme on an interval $[t_0,t_1]$ for $t_1>t_0$. First define
\begin{equation}
v^{\rho,r,t_0,(0)}_i:=v^{\rho,r,\nu}_i(t_0,.)\ast_{sp}G^{\rho,r}_{\nu}.
\end{equation}
Note that we have
\begin{equation}
\begin{array}{ll}
v^{\rho,r,\nu,t_0,(0)}_i=v^{\rho,r,\nu}_i(t_0,.)\ast_{sp}G^{\rho,r}_{\nu}\\
\\
=\left( h_i\ast_{sp}G^{\rho,r}_{\nu}(t_0,.)\right) 
G^{\rho,r}_{\nu}(t_1-t_0,.)+
{\big (}\sum_{k\geq 1}\delta v^{\rho,r,\nu,(k)}_i{\big )}\ast_{sp}G^{\rho,r}_{\nu}(t_1-t_0,.)\\
\\
=h_i\ast_{sp}G^{\rho,r}_{\nu}(t_1,.)+{\big (}\sum_{k\geq 1}\delta v^{\rho,r,\nu,(k)}_i(t_0,.){\big )}\ast_{sp}G^{\rho,r}_{\nu}(t_1-t_0)\in {\cal C}^{m(n+1)-\mu}_{pol,m},
\end{array}
\end{equation} 
where for the second summand we may use (\ref{assdelta}).
For $k\geq 1$ define on $[t_0,t_1]$
\begin{equation}\label{Navleraysol}
\begin{array}{ll}
v^{\rho,r,\nu,t_0,(k)}_i
=v^{\rho,r,\nu}(t_0,.)\ast_{sp}G^{\rho,r}_{\nu}+\sum_{j=1}^n \left( v^{\rho,r,\nu,t_0,(k-1)}_j\frac{\partial v^{\rho,r,\nu,t_0,(k-1)}_i}{\partial x_j}\right) \ast G_{\nu}\\
\\ \hspace{1cm}+\sum_{j,m=1}^n\int_{{\mathbb R}^n}\left( \frac{\partial}{\partial x_i}K_n(.-y)\right) \sum_{j,m=1}^n\left( \frac{\partial v^{t_0,\rho,r,\nu,(k-1)}_m}{\partial x_j}\frac{\partial v^{t_0,\rho,r,\nu,(k-1)}_j}{\partial x_m}\right) (.,y)dy\ast G_{\nu}.
\end{array}
\end{equation}
 Define
\begin{equation}
\delta v^{t_0,\rho,r,\nu,(0)}_i(t,x):=0,
\end{equation}
and for $k\geq 1$
\begin{equation}\label{Navleraysol}
\begin{array}{ll}
\delta v^{t_0,\rho,r,\nu,(k)}_i:=v^{t_0,\rho,r,\nu, (k)}_i-v^{t_0,\rho,r,\nu,(k-1)}_i.
\end{array}
\end{equation}
For $\tau \geq 0$ and $m\geq 2$ we get
\begin{equation}
\delta v^{\rho,r,\nu,t_0,(k)}(\tau,.)\in {\cal C}^{m(n+1)}_{pol,m},~v^{\rho,r,\nu,t_0,(k)}(\tau,.)\in {\cal C}^{m(n+1)-\mu}_{pol,m}
\end{equation}
for some $\mu\in (0,1)$. This observation can be used to designing global schemes using damping effects of Gaussian convolutions auto-control, or external control.
\end{rem}

Next we observe that we can strengthen the observations above if we avoid the convolution of the initial data in the local scheme. We go back to the unscaled version of the scheme for simplicity of notation (i.e. $\rho=r=1$ in the following).
Assume at time $t_0\geq 0$ we have Cauchy data 
\begin{equation}
v^{\nu,0}_i(t_0,.):=v^{\nu}_i(t_0,.)\in {\cal C}^{m(D+1)}_{pol,m}
\end{equation}
For $0\leq |\beta|\leq m$ and $|\gamma|+1=|\beta|,~\gamma_{j}+1=\beta_j$ if $|\beta|>0$ define the local time iteration scheme
\begin{equation}\label{Navlerayschemeiii}
\begin{array}{ll}
 D^{\beta}_xv^{\nu,k}_i=D^{\beta}_xv^{\nu}_i(t_0,.)
-D^{\gamma}_x\left( \sum_{j=1}^D v^{\nu,k-1}_j\frac{\partial v^{\nu,-,k-1}_i}{\partial x_j}\right) \ast G_{\nu,j}\\
\\+\left( \sum_{j,m=1}^D\int_{{\mathbb R}^D}D^{\gamma}_x\left( \frac{\partial}{\partial x_i}K_D(.-y)\right) \sum_{l,m=1}^D\left( \frac{\partial v^{\nu,k-1}_m}{\partial x_l}\frac{\partial v^{\nu,k-1}_l}{\partial x_m}\right) (t,y)dy\right) \ast G_{\nu,j}.
\end{array}
\end{equation}
Here, recall $G_{\nu}$ is the fundamental solution of the heat equation $p_{,t}-\nu\Delta p=0$, $\ast$ denotes the convolution, $\ast_{sp}$ denotes the spatial convolution, and $K_D$ denotes the fundamental solution of the 
Laplacian equation for dimension $D\geq 3$. In the following the constant $c>0$ is generic. For $1\leq i\leq D$ the initial data $v^{\nu}_i(t,.)$ are in ${\cal C}^{m(D+1)}_{pol,m}$. Hence, for $k=0$, for $0\leq |\gamma|\leq m$, and for $|x|\geq 1$
\begin{equation}
{\big |}D^{\gamma}_xv^{\nu,0}_i(t_0,x){\big |}\leq \frac{c}{1+|x|^{m(D+1)}}
\end{equation}
 for some finite constant $c>0$ and $t_0\geq 0$. 
Assuming inductively that for $t\in [t_0,t_0+\Delta]$
\begin{equation}\label{indhyp}
\forall~l\leq k-1~\forall 0\leq|\gamma|\leq m~{\Big |}D^{\gamma}_xv^{\nu,l}_i(t,.){\Big |}\leq  \frac{c}{1+|x|^{m(D+1)}}
\end{equation}
we have or some finite constant $c>0$, for $0\leq |\delta|\leq m-1$  and for $|x|\geq 1$ 
\begin{equation}
{\big |}D^{\delta}_xB^{k-1}{\big |}:={\Big |}\sum_{j=1}^D D^{\delta}_x\left( v^{\nu,k-1}_j\frac{\partial v^{\nu,k-1}_i}{\partial x_j}\right) (t,.){\Big |}\leq \frac{c}{1+|x|^{2m(D+1)}},
\end{equation}
and 
\begin{equation}
{\big |}D^{\delta}_xL^{k-1}{\big |}\leq \frac{c}{1+|x|^{2m(D+1)-1}},
\end{equation}
where 
\begin{equation}
D^{\delta}_xL^{k-1}\equiv\sum_{j,m=1}^D\int_{{\mathbb R}^D}\left( \frac{\partial}{\partial x_i}K_D(.-y)\right) \sum_{j,m=1}^D\left( D^{\delta}_x\left( \frac{\partial v^{\nu,k-1}_m}{\partial x_j}\frac{\partial v^{\nu,k-1}_j}{\partial x_m}\right)\right)  (t,y)dy.
\end{equation}
Convolutions with  $G_{\nu}$ or $G_{\nu,i}$ weaken this polynomial decay by order $D$ at most such that we (generously) get for some finite constant $c>0$, for $0\leq |\delta|\leq m-1$  and for $|x|\geq 1$ 
\begin{equation}\label{B}
{\big |}D^{\delta}_xB^{k-1}\ast G_{\nu,j}{\big |}\leq \frac{c}{1+|x|^{(2m-1)(D+1)}}
\end{equation}
and
\begin{equation}\label{L}
{\Big |}D^{\delta}_xL^{k-1}\ast G_{\nu,j}{\Big |}\leq \frac{c}{1+|x|^{(2m-1)(D+1)-1}}.
\end{equation}
Hence using the representation (\ref{Navlerayscheme}),(\ref{B}),(\ref{L}) we can complete the  induction step and get
\begin{equation}\label{ind}
\forall~l\leq k~\forall 0\leq|\gamma|\leq m~{\big |}D^{\gamma}_xv^{\nu,l}_i(.){\big |}\leq  \frac{c}{1+|x|^{m(D+1)}}
\end{equation}
and by (\ref{indhyp}) the same holds for the increments $$D^{\gamma}_x\delta v^{\nu,k}_i(.)=D^{\gamma}_xv^{\nu,k}_i(.)-D^{\gamma}_x v^{\nu,k-1}_i(.).$$
For some $\Delta >0$ we have local time contraction with respect to a $H^m\cap C^m$-norm, such that the limit
\begin{equation}
D^{\gamma}_x v^{\nu}_i=\lim_{k\geq 1}D^{\gamma}_x v^{\nu,k}_i
\end{equation}
inherits this order of polynomial decay at spatial infinity for all $0\leq |\gamma|\leq m$. More precisely, for $t\in [t_0,t_0+\Delta]$, $|x|\geq 1$, and for all $0\leq |\gamma|\leq m$ we have a finite constant $c>0$ such that for all $x\in {\mathbb R}^D$
\begin{equation}\label{polybound}
\forall 0\leq|\gamma|\leq m~{\Big |}D^{\gamma}_x v^{\nu}_i(t,x){\Big |}\leq  \frac{c}{1+|x|^{m(D+1)}}.
\end{equation} 
% and
% \begin{equation}\label{polybound2}
% \forall~ k\geq 1~\forall 0\leq|\gamma|\leq 3~{\Big |}D^{\gamma}_x\delta v^{\nu,-,k}_i(t,x){\Big |}\leq  \frac{c}{1+|x|^{2(D+1)}}.
% \end{equation} 
%  For small $T>0$ it is observed that the constant $c>0$ can be chosen independently of the iteration index $k$. Furthermore, for $0\leq |\gamma|\leq 4$ 
% \begin{equation}
% {\Big |}v^{f}_i\ast_{sp} D^{\gamma}_xG_{\nu}(t,.){\Big |}\leq  \frac{c}{1+|x|^{2(D+1)}}~\mbox{$t>0$ and $c$ depends on $t>0$},
% \end{equation}
% such that a corresponding statement holds for the value functions $v^{\nu,-,k}_i,~1\leq i\leq 3$ and $k\geq 2$.
% Given $k\geq 1$ choose a number $m$ is such that for $t>0$
% \begin{equation}
% \forall~0\leq |\gamma|\leq m~D^{\gamma}_x\delta v^{init,\nu,-,k}_i(t,.) \mbox{is continuous and bounded.}
% \end{equation}
% For $0\leq |\gamma|\leq m$ and for $|\beta|+1=|\gamma|,~\beta_j+1=\gamma_j$ (if $|\gamma|\geq 1$) we consider a representation of  $D^{\gamma}_xv^{\nu,-,k}_i,~1\leq i\leq D,~k\geq 1$ of the form

We choose a sequence $(\nu_p)_{p\geq 1}$ converging to zero and consider the spatial transformation
\begin{equation}
v^{c,\nu_p}_i(t,y)=v^{\nu_p}_i(t,x)
\end{equation}
for $y_j=\arctan(x_j),~1\leq j\leq D$ and for all $t\in [t_0,t_0+\Delta]$. 
For multiindices $\gamma$ with $0\leq |\gamma|\leq m$ for all $t\in [0,T]$, all $y\in \left(-2\pi ,2\pi \right)^D$ and all $x\in {\mathbb R}^D$ 
\begin{equation}
{\big |}D^{\gamma}_y\delta v^{c,\nu_p}_i(t,y){\big |}\leq c_0(1+|x|^{2m}){\big |}D^{\gamma}_x\delta v^{\nu_p}_i(t,x){\big |}\leq C
\end{equation}
for some finite constants $c_0,C>0$. Then for and $\epsilon >0$ there is a subsequence which we may denote again by $(\nu_p)_{p\geq 1}$ such that we have a limit
\begin{equation}
\lim_{p\uparrow \infty}v^{c,\nu_p}_i(t,.):=\lim_{\nu_p\downarrow 0} v^{c,\nu_p}_i(t,.)\in H^{m-\epsilon}~\mbox{ for all $t_0\leq t\leq t_0+\Delta$}
\end{equation}
by Rellich's theorem.
Moreover $v^{\nu_p}_i(t,.)\in C^{m-1}_0\left(\left(-\frac{\pi}{2},\frac{\pi}{2}\right)^D  \right) $. Here, $$C^{m-1}_0\left(\left(-\frac{\pi}{2},\frac{\pi}{2}\right)^D  \right)$$ is the function space of $m-1$-times continuously differentiable functions which vanish if a component $y_i$ becomes equal to $-\frac{\pi}{2}$ or $\frac{\pi}{2}$. This space is a Banach space if equipped with the usual $C^{m-1}$-supremum-norm on the bounded domain. 
The limit $e^c_i(t,.)=\lim_{p\uparrow \infty}v^{c,\nu_p}_i(t,.)\in H^{m-\epsilon}\cap C^{m-1}, 1\leq i\leq D,~t\in [t_0,t_0+\Delta_0]$ satisfies the transformed Euler equation with respect to spatial coordinates $y_i,~1\leq i\leq D$, and the corresponding function limit $e_i, 1\leq i\leq D$ with $e_i(t,x)=e^{c}_i(t,y)$ satisfies the original Euler equation with respect to spatial coordinates $x_i,~1\leq i\leq D$. The construction becomes global straightforwardly by application of the semigroup property.

\section{Short and Long time singularities of the Euler equation}
A characteristic difference of the Euler equation (compared to the Navier Stokes equation) is that it can be solved backwards in time for regular data (at least locally). 
In the previous section we have observed that compactness arguments can be applied if we have strong polynomial decay at spatial infinity of the data. This strong polynomial decay is inherited by the natural local time iteration scheme considered above. We have observed that local-time contraction results hold in the viscosity limit for regular data. In the convolution representation of the velocity component functions the contribution of the Gaussian or its first order spatial derivative is concentrated on a ball of radius $\sqrt{\nu}$ around the spatial argument $x$ for positive viscosity $\nu >0$ (in original coordinates).  We have observed that subsequences in strong spaces with uniform upper bounds (independent of $\nu$) have natural pointwise limits. Note that spatial Fourier transforms of the convolution of regular data with (first order spatial derivatives) Gaussian have the effect of a multiplication with a linear term in the viscosity limit (at most) which is absorbed by the polynomial spatial decay of the Fourier transform of the regular data. If short time solutions $e^{-}_i,~1\leq i\leq D$ of the time-reversed Euler equation with weakly singular data gain 'enough' regularity then this implies the existence of weak short - and even long-time singularities of the original Euler equation. Here 'enough regularity' for long-time singularities means that the evaluation of a local time solution with weakly singular data at some time has sufficient regularity such that the semi-group property of the (time-reversed) Euler equation can be combined  with the argument for a global solution branch of the previous section. Next we consider this in detail. We consider positive viscosity $\nu >0$ first, and consider the viscosity limit in a second step. We consider the time transformation $t\rightarrow -t=:s$ and the time reversed equation for
\begin{equation}
v^{\nu,-}_i(s,.)=~1\leq i\leq D,~v^{\nu,-}_i(s_0,.)=v^{\nu}_i(t_1,.)~1\leq i\leq D.
\end{equation}
for some $t_1>0$, where 
\begin{equation}\label{Navlerayequationr10}
\begin{array}{ll}
 \frac{\partial v^{\nu -}_i}{\partial t}-\nu \Delta v^{\nu}_i
-\sum_{j=1}^D \left( v^{\nu,-}_j\frac{\partial v^{\nu,-}_i}{\partial x_j}\right) \\
\\+\int_{{\mathbb R}^D} K_{D,i}(.-y)\sum_{j,m=1}^D\left( \frac{\partial v^{\nu,-}_m}{\partial x_j}\frac{\partial v^{\nu,-}_j}{\partial x_m}\right) (.,y)dy=0.
\end{array}
\end{equation}
 
Since we have global solution branches for the class of strong data ${\cal C}^{m(D+1)}_{pol,m}$ (cf. previous section) 
the initial time $s_0=-t_1$ for the time-reversed equation can be chosen arbitrarily if a short time solution $v^{\nu,-}_i,~1\leq i\leq D$ with weakly singular data $v^{\nu,-}_i(s_0,.),~1\leq i\leq D$ gains enough regularity after short time such that $v^{\nu,-}_i(s,.)\in {\cal C}^{m(D+1)}_{pol,m}$ for some $s>s_0$. We reconsider here a variation of a local construction which we have considered elsewhere and sharpen some results.
We construct local time solutions for carefully chosen data via the iteration scheme $v^{\nu,-,k}_i,~1\leq i\leq D,~k\geq 0$, where
\begin{equation}\label{Navlerayschemekk}
\begin{array}{ll}
 v^{\nu,-,k}_i=v^{\nu,-}_i(s_0,.)\ast_{sp}G_{\nu}
+\sum_{j=1}^D \left( v^{\nu,-,k-1}_j\frac{\partial v^{\nu,-,k-1}_i}{\partial x_j}\right) \ast G_{\nu}\\
\\-\left( \sum_{j,m=1}^D\int_{{\mathbb R}^D}\left( \frac{\partial}{\partial x_i}K_D(.-y)\right) \sum_{j,m=1}^D\left( \frac{\partial v^{\nu,-,k-1}_m}{\partial x_j}\frac{\partial v^{\nu,-,k-1}_j}{\partial x_m}\right) (.,y)dy\right) \ast G_{\nu},~k\geq 2\\
\\
v^{\nu,-,0}_i(s_0,.):=v^{\nu,-}_i(s_0,.):=v^{\nu}_i(t_1,.),~1\leq i\leq D,\\
\\
v^{\nu,-,1}_i(s_0,.):=v^{\nu,-}_i(s_0,.)\ast_{sp}G_{\nu},~1\leq i\leq D.
\end{array}
\end{equation}
Local time contraction results are obtained as for the original Navier Stokes equation, and we may use representations of solutions of the form
\begin{equation}
v^{\nu,-}_i=v^{\nu,-,2}_i(s_0,.)+\sum_{k\geq 3}\delta v^{\nu,-,k}_i,
\end{equation}
where we denote $\delta v^{\nu,-,k}_i=v^{\nu,-,k}_i-v^{\nu,-,k-1}_i$ for $k\geq 1$.
We choose data which are weakly singular in the sense that the vorticity (original Euler equation)
\begin{equation}
\omega =\mbox{curl}(v)=\left(\frac{\partial v_3}{\partial x_2}-\frac{\partial v_2}{\partial x_3},\frac{\partial v_1}{\partial x_3}-\frac{\partial v_3}{\partial x_1},\frac{\partial v_2}{\partial x_1}-\frac{\partial v_1}{\partial x_2} \right)
\end{equation}
has no finite upper bound at time $t_1>0$. This means that there are data $\mbox{curl}(h)\in {\cal C}^{m(D+1}_{pol,m}$  such that a solution of the Cauchy problem for the incompressible Euler equation in vorticity form
\begin{equation}\label{vorticity}
\frac{\partial \omega}{\partial \tau}+v\cdot \nabla \omega=\frac{1}{2}\left(\nabla v+\nabla v^T\right)\omega, 
\end{equation}
blows up after finite time $t_1$, where $t_1>0$ can be large. Note that a vorticity blow up means that the corresponding velocity solution has a kink as it is well-known (cf. \cite{MB} )that
\begin{equation}\label{vel}
v(t,x)=\int_{{\mathbb R}^3}K_3(x-y)\omega(t,y)dy,~\mbox{where}~K_3(x)h=\frac{1}{4\pi}\frac{x\times h}{|x|^3}.
\end{equation}

We prove
\begin{thm}\label{main1}
Let $D= 3$. For any (arbitrarily large) finite time $T>0$ and all $m\geq 2(D+1)$ there exist data $h_i\in {\cal C}^{m(D+1)}_{pol,m},~1\leq i\leq D$ and a vorticity solution $\omega_i,~1\leq i\leq D$ of the $D$-dimensional incompressible Euler equation Cauchy problem such that there is a blow-up of the classical solution at time $T>0$, i.e.,
\begin{itemize}
 \item[i)] there is a solution function $\omega_i:[0,T)\times {\mathbb R}^D\rightarrow {\mathbb R},~1\leq i\leq D$ in $C^{1}\left(\left[0,T\right), {\cal C}^m_{\mbox{pol}} \right)$ which satisfies the incompressible Euler equation pointwise on the domain $\left[0,T\right)\times {\mathbb R}^D$ in a classical sense;
 \item[ii)] for the solution in item i) we have
 \begin{equation}\label{supomega}
 \sup_{t\in [0,T)}|\omega_i(t,x)|= \infty,
 \end{equation}
i.e., there is no finite upper bound for the left side of (\ref{supomega}). 
\end{itemize}
\end{thm}
% \begin{defi}
% We say that a classical solution branch $\omega_i,~1\leq i\leq D$ of an incompressible Euler Cauchy problem with data $\omega^f_{i},~1\leq i\leq D$ with  $\omega^f_i\in C^{m}$ for $m\geq k$ has spatial a kink of order $k\geq 1$, if there is a space-time point $(\tau,x)$ with $\tau>0$ and $x\in {\mathbb R}^D$ such that
% \begin{equation}
% \omega_i\in C^{k}\setminus C^{k-1}~\mbox{at }~(\tau,x)
% \end{equation}
% for some integer $k\geq 1$.
% \end{defi}
% 
% 
% \begin{cor}\label{maincor}
% Let $D=3$. For any $k\geq 2$ and $s\geq 0$ there exists data $h_i\in H^m\cap C^{m}$ with $m\geq k+2$ and a vorticity solution $\omega_i,~1\leq i\leq D$ of the three dimensional incompressible Euler equation Cauchy problem with data $h_i,~1\leq i\leq n$ such that after some finite time the solution has a kink of order $k$.
% \end{cor}
Note that the preceding theorem and the construction of global solution branches for data $h_i\in {\cal C}^{m(D+1)}_{pol,m},~1\leq i\leq D,~m\geq 2$ in the preceding section imply that incompressible Euler equation Cauchy problems do not have unique solutions in general.
More precisely we have
\begin{cor}
For data $D\geq 3$ $h_i\in {\cal C}^{m(D+1)}_{pol,m},~1\leq i\leq D,~m\geq 2$ the Cauchy problem of the incompressible Euler equation as infinitely many solutions. Next to a global solution branch there exist solution with blow ups of first derivatives at any time $T>0$ and solution in $C^{k+1}\setminus C^k, k\geq 2$ or solutions with kinks of any order $k$ at any time $T>0$.  
\end{cor}
Some arguments of the preceding section such as local time contraction can be transferred to the time-reversed equation straightforwardly. Additionally we have to show that for some weakly singular data in $H^2$ there is a local solution branch which gains enough regularity after short time in order to apply the arguments for global regular solution branches of the Euler equation for regular data obtained in the last section. 
Finally, we add four additional steps which we need for a proof of Theorem \ref{main1}.
\begin{itemize}
 \item[i)] First we choose appropriate weakly singular data. For some time $s_0$, a positive real number $\beta_0\in (1,1+\alpha_0)$, and $\alpha_0\in \left(0,\frac{1}{2}\right)$ we consider  velocity component data $v^{\nu,-}_i(s_0,.)\in H^{2},~1\leq i\leq 3$. For one index $i_0\in \left\lbrace 1,2,3\right\rbrace$  we  choose weak data, where we define $v^{\nu,-}_{i_0}(s_0,x)=g_{(0)}(r)$ for some univariate function $g$. The function $g_{(0)}$ depends on $r=\sqrt{x_1^2+x_2^2+x_3^2}\geq 0$. We define  $g_{(0)}:{\mathbb R}^0_+\rightarrow {\mathbb R}$ by
\begin{equation}
g_{(0)}(r):= \left\lbrace \begin{array}{ll}\phi_{1}(r)r^{\beta_0}\sin\left(\frac{1}{r^{\alpha_0}}\right)\\
0
              \
             \end{array}\right.
\end{equation}
where $\phi \in C^{\infty}_c$, and
\begin{equation}
\phi_{1}(r)=\left\lbrace \begin{array}{ll}
1~\hspace{1.75cm}\mbox{if }~r\leq 1,\\
\phi_1(r)=\alpha_*(r) \mbox{ if }~1\leq r\leq 2,
\\
0~\hspace{1.65cm}\mbox{ if }~r\geq 2.
\end{array}\right.
\end{equation}
Here, $\alpha_*$ is a smooth function with bounded derivatives for $1\leq r\leq 2$, and $C^{\infty}_c$ denotes the function space of smooth functions with compact support. 
For $j\in \left\lbrace 1,2,3\right\rbrace\setminus \left\lbrace i_0 \right\rbrace$ we may choose regular velocity component data, i.e. we choose data 
\begin{equation}
v^{\nu,-}_j(s_0,.)\in C^{\infty}_{c0}.
\end{equation}
Note that for all $1\leq i\leq D$ and multiindices $\alpha$ with $0\leq |\alpha|=k$ we have
\begin{equation}
{\big |}D^{\alpha}_xv^{\nu,-}_i(s_0,.){\big |}\leq C r^{\beta_0-k(1+\alpha_0)}.
\end{equation}
For the derivative of the data $v^{\nu,-}_{i_0}(s_0,.)$ we compute for $r\neq 0$ and $r\leq 1$
\begin{equation}
\begin{array}{ll}
g'(r)=\frac{d}{dr} r^{\beta_0}\sin\left(\frac{1}{r^{\alpha_0}}\right)
=\beta_0 r^{\beta_0-1}\sin\left(\frac{1}{r^{\alpha_0}}\right)-\alpha_0 r^{\beta_0-1-\alpha_0}\cos\left(\frac{1}{r^{\alpha_0}}\right).
\end{array}
\end{equation}
The derivative $g'$ of the function $g$ at $r=0$ is strongly singular for $\beta\in (1,1+\alpha_0)$ and $\alpha_0\in \left(0,\frac{1}{2}\right)$. Note that it is 'oscillatory' singular bounded for $\beta_0 -1=\alpha_0 \in \left(0,\frac{1}{2}\right)$. Note that for data $v^{\nu,-}_{i_0}(s_0,x_1,x_2,x_3)=g(r)$ we have (for $r\neq 0$)
\begin{equation}
v^{\nu,-}_{i_0,j}(s_0,x)=g'(r)\frac{\partial r}{\partial x_j}=g'(r)\frac{x_j}{r}.
\end{equation}
In polar coordinates $(r,\theta,\phi)\in [0,\infty)\times [0,\pi]\times[0,2\pi]$ with 
\begin{equation}
x_1=r\sin(\theta)\cos(\phi),~x_2=r\sin(\theta)\sin(\phi),~x_3=r\cos(\theta),
\end{equation}
(where for $r\neq 0$ and $x_1\neq 0$ we have $r=\sqrt{x_1^2+x_2^2+x_3^2}, ~\theta=\arccos\left(\frac{x_3}{r}\right),~\phi=\arctan\left(\frac{x_2}{x_1} \right)$) we get
\begin{equation}
\begin{array}{ll}
v^{\nu,-}_{i_0,1}(s_0,x)=g'(r)\frac{x_1}{r}=g'(r)\sin(\theta)\cos(\phi),\\
v^{\nu,-}_{i_0,2}(s_0,x)=g'(r)\frac{x_2}{r}=g'(r)\sin(\theta)\sin(\phi),\\
v^{\nu,-}_{i_0,3}(s_0,x)=g'(r)\frac{x_3}{r}=g'(r)\cos(\theta),
\end{array}
\end{equation}
such that we have
\begin{equation}
v^{\nu,-}_{i_0}(s_0,.)\in H^{1} ~\mbox{obviously.}
\end{equation}
The second derivative of $g$ is
\begin{equation}
\begin{array}{ll}
g''(r)=\frac{d^2}{dr^2} r^{\beta_0}\sin\left(\frac{1}{r^{\alpha_0}}\right)\\
\\
=\frac{d}{dr}\left( \beta_0 r^{\beta_0-1}\sin\left(\frac{1}{r^{\alpha_0}}\right)-\alpha_0 r^{\beta_0-1-\alpha_0}\cos\left(\frac{1}{r^{\alpha_0}}\right)\right)\\
\\
=\beta_0(\beta_0-1)r^{\beta_0-2}\sin\left(\frac{1}{r^{\alpha_0}}\right)
-\alpha_0\beta_0 r^{\beta_0-3-\alpha_0}\cos\left(\frac{1}{r^{\alpha_0}}\right)\\
\\
+(\alpha_0)(1+\alpha_0-\beta_0) r^{\beta_0-1-\alpha_0}\cos\left(\frac{1}{r^{\alpha_0}}\right)\\
\\
-(\alpha_0)^2 r^{\beta_0-2-2\alpha_0}\sin\left(\frac{1}{r^{\alpha_0}}\right).
\end{array}
\end{equation}
We have $v^{\nu,-}_{i_0}(s_0,.)\in H^{2}$, since
\begin{equation}
\beta_0-2-2\alpha_0 >-\frac{3}{2}.
\end{equation}
Note that
\begin{equation}
v^{\nu,-}_{i_0}(s_0,.) \in C^{\delta}\left({\mathbb R}^3\right), 
\end{equation}
for H\"{o}lder constants of order $\delta \in \left(0,\beta-\alpha\right)$. Hence Lipschitz continuity of the data does not hold, and we have to refine estimates of convolutions with (first spatial derivatives of) the Gaussian in order to extend the argument of the preceding section.

\item[ii)] For the scheme defined in (\ref{Navlerayschemekk}) above and data
$v^{\nu,-,0}_i(s_0,.),~1 \leq i\leq D$ as defined in item i) we first observe
(based on similar reasons as in the previous section) that for first order multivariate spatial derivatives $0\leq |\gamma |\leq 1$ and $|x|\geq 1$, for some $\Delta >0$ for all $s \in [s_0,s_0+\Delta]$ there exists a finite constant $c$ which is independent of $\nu >0$ such that
\begin{equation}\label{polybound}
~~{\Big |}D^{\gamma}_x v^{\nu,-,2}_i(s,x){\Big |}\leq  \frac{c}{1+|x|^{2(D+1)}},
\end{equation} 
and
\begin{equation}\label{polybound2}
\forall~ k\geq 1~{\Big |}D^{\gamma}_x v^{\nu,-,k}_i(s,x){\Big |}\leq  \frac{c}{1+|x|^{2(D+1)}}.
\end{equation} 
Note that this implies that for first order multivariate spatial derivatives $0\leq |\gamma |\leq 1$, some $\Delta >0$ and $s\in [s_0,s_0+\Delta]$ and $|x|\geq 1$ there exists a finite constant $c$ which is independent of $\nu >0$ such that for $\delta v^{\mbox{init},\nu,-,2}_i:=v^{\nu,-,2}_i-v^{\nu,-,2}_i(s_0,.)\ast_*G_{\nu}$ we have
\begin{equation}\label{polybounddelta}
~~{\Big |}D^{\gamma}_x\delta v^{\mbox{init},\nu,-,2}_i(s,x){\Big |}\leq  \frac{c}{1+|x|^{2(D+1)}},
\end{equation} 
and
\begin{equation}\label{polybound2delta}
\forall~ k\geq 1~~{\Big |}D^{\gamma}_x\delta v^{\nu,-,k}_i(s,x){\Big |}\leq  \frac{c}{1+|x|^{2(D+1)}}.
\end{equation} 
Next concerning the behavior at $r=\sqrt{\sum_{i=1}^Dx_i^2}=0$ we have to refine some observations of the preceding section since the data are only H\"{o}lder continuous with exponent $\delta\in (0,\beta_0-\alpha_0)$ and not Lipschitz in general as $\beta_0-\alpha_0$ is close but smaller than $1$ according to our choice. A choice $\beta_0=1+\alpha_0$ implies a bounded oscillatory singularity for the vorticity. Here we want to prove a stronger result, i.e., a blow up for vorticity. We consider convolutions of (first order spatial derivatives of) the Gaussian $G_{\nu}$ with H\"{o}lder continuous functions $g$, where we have the initial data in mind first. It seems impossible to get $v$-independent estimates for the increments of the first iteration of the local scheme, but the local functional increments have $\nu$-independent estimates from the second iteration step on, and this insufficient for our purposes. Indeed, for some $\Delta>0$ and for all $s\in [s_0,s_0+\Delta]$ we have $\nu$-independent upper bounds
\begin{equation}\label{polybounddelta0}
~~{\Big |}D^{\gamma}_x\delta v^{\mbox{init},\nu,-,2}_i(s,x){\Big |}\leq  r^{\beta_0-|\gamma|},
\end{equation} 
and
\begin{equation}\label{polybound2delta0}
\forall~ k\geq 3~~{\Big |}D^{\gamma}_x\delta v^{\nu,-,k}_i(s,x){\Big |}\leq  r^{\beta_0-|\gamma|}.
\end{equation} 
From these estimates in (\ref{polybounddelta0}), (\ref{polybound2delta0}), (\ref{polybounddelta}), and (\ref{polybound2delta}) local time contraction can be obtained straightforwardly. Then it follows that
\begin{equation}
v^{\nu,-}_i=v^{\nu,-}_i(s_0,.)+\delta v^{\mbox{init},\nu,-,2}_i(s,.)+\sum_{k=3}^{\infty}\delta v^{\nu,-,k}_i(s,x)\in C^{1,2}
\end{equation}
has a uniform upper bound such that
\begin{equation} 
e^{-}_{i}=\lim_{k\uparrow \infty }v^{\nu_k,-}_i\in C^{1,1}\cap {\cal C}^{m(D+1)}_{pol,1}.
\end{equation}
Let is consider this in some more detail.  Since we consider local time solutions we may consider representations of local solutions in terms of scaled Gaussian $G^{\rho,r}_{\nu}\equiv G_{\nu'}$ as above with $\nu'=4\pi \rho r^2 \nu$. Recall that we start the iteration with the initial data
\begin{equation}
v^{\nu',-,0}_i(s_0,.)=v^{\nu',-}_i(s_0,.),~1\leq i\leq D
\end{equation}
assuming that the data $v^{\nu',-}_i(s_0,.),~1\leq i\leq D$ are known at time $s_0\geq 0$, and that for $k=1$ we define
\begin{equation}
v^{\nu',-,1}_i=v^{\nu',-}_i(s_0,.)\ast G_{\nu'},~1\leq i\leq D
\end{equation}
In order to estimate $\delta v^{\nu',-,2}_i=\delta v^{\mbox{init},\nu',-,2}_i= v^{\nu',-,2}_i- v^{\nu',-,0}_i(s_0,.)\ast_{sp}G_{\nu'}$ and $\delta v^{\nu',-,k}_i,\geq 3$ for $1\leq i\leq D$ we have to plug in 

$v^{\nu',-,1}_i=v^{\nu',-}_i(s_0,.)\ast G_{\nu'},~1\leq i\leq D$ into the Burgers term and the Leray projection term. 
 Using the convolution rule we observe that for $s_1>s_0$ the first order spatial derivatives of
\begin{equation}
\left( v^{\nu,-}_i\ast G_{\nu'}\right) (s_1,x)=\int_{s_0}^{s_1}\int_{{\mathbb R}^D} v^{\nu',-}(s_0,x-y)\frac{1}{\sqrt{4\pi\nu' \sigma}^D}
\exp\left(-\frac{|y|^2}{4\nu' \sigma} \right)dyd\sigma
\end{equation}
have the representation
\begin{equation}
\left( v^{\nu,-}_i(s_0,.)\ast G_{\nu' ,i}\right) (s_1,x)=\int_{s_0}^{s_1}\int_{{\mathbb R}^D} v^{\nu',-}_i(s_0,x-y)G_{\nu,i}(\sigma,y)dyd \sigma.
\end{equation}
Now the first order derivatives of the scaled Gaussian $G_{\nu'}$ are given by
\begin{equation}
G_{\nu',i}(s,x)=\left( \frac{-2y}{4\rho r^2\nu \sigma}\right) \frac{1}{\sqrt{4\pi\nu' \sigma}^D}\exp\left(-\frac{-|y|^2}{4\nu' \sigma} \right),
\end{equation}
and have the upper bound
\begin{equation}
\begin{array}{ll}
{\big |}G_{\nu',i}(\sigma,y){\big |}\leq {\Big |}\left( \frac{-2y}{4\nu' \sigma}\right) \frac{1}{\sqrt{4\pi\nu' \sigma}^D}\exp\left(-\frac{-|y|^2}{4\nu' \sigma} \right){\Big |}\\
\\
\leq {\Big |}\frac{2}{(4\pi\nu'\sigma)^{\delta}}
\frac{1}{|y|^{D+1-2\delta}} \left( \frac{|y|^2}{4\pi\nu' \sigma}\right)^{D/2+1-\delta}\exp\left(-\frac{-|y|^2}{4\nu' \sigma} \right) {\Big |}\\
\\
\leq {\Big |}\frac{1}{(4\pi\nu'\sigma)^{\delta}}
\frac{c_s}{|y|^{D+1-2\delta}}{\Big |},
\end{array}
\end{equation}
where for $|z|=\frac{|y|}{\sqrt{4\pi\nu' \sigma}}$ and $\delta \in \left(0,1\right)$ 
\begin{equation}
c_s:=2\sup_{|z|\geq 0}\left( \frac{|y|^2}{4\pi\rho r^2\nu \sigma}\right)^{D/2+1-\delta}\exp\left(\frac{-|y|^2}{4\nu' \sigma} \right).
\end{equation}
Since
\begin{equation}
{\big |}v^{\nu',-}(s_0,x-y){\big |}\leq c_0|x-y|^{\beta_0}
\end{equation}
for some finite constant $c_0>0$ we have 
\begin{equation}\label{uppb}
\begin{array}{ll}
{\big |}v^{\nu',-}(s_0,.)\ast G_{\nu ,i}{\big |}\leq 
c_1|r|^{\beta_0-1}
\end{array}
\end{equation}
for some finite constant $c_1$.
Here, note that the elliptic integral gives the upper bound $\int_{s_0}^{s_1}\frac{r^{\beta_0+2\delta-1}}{(\nu'(\sigma-s_0))^{\delta}}d\sigma$, and as we plug in this elliptic integral into the iteration formula for $v^{\nu',2,-}_i$ the contribution the of the integral terms in the formula for $v^{\nu',2,-}_i$ is concentrated in the area $r^2\leq \nu'$ as $\nu'$ becomes small. Concerning the behavior of the first order spatial derivatives of the Gaussian for $\sqrt{\sum_{i=1}^Dy_i^2}=r>\sqrt{\nu'}$ note that 
\begin{equation}\label{firstderg}
\begin{array}{ll}
{\big |}G_{\nu',i}(\sigma,y){\big |}
\leq {\Big |}\left( \frac{-2y}{4\nu' \sigma}\right) \frac{1}{\sqrt{4\pi\nu' \sigma}^D}\exp\left(-\frac{-|y|^2}{8\nu' \sigma} \right)\exp\left(-\frac{-|y|^2}{8\nu' \sigma} \right){\Big |}\\
\\
\leq {\Big |}\frac{1}{(4\pi\nu'\sigma)^{\delta}}
\frac{\tilde{c}_s}{|y|^{D+1-2\delta}}\exp\left(-\frac{-|y|^2}{8\nu' \sigma} \right){\Big |},
\end{array}
\end{equation}
where 
\begin{equation}
\tilde{c}_s:=2\sup_{|z|\geq 0}\left( \frac{|y|^2}{4\pi\rho r^2\nu \sigma}\right)^{D/2+1-\delta}\exp\left(\frac{-|y|^2}{8\nu' \sigma} \right).
\end{equation} 
Convolutions of this upper bound with data of order $|x-y|^{\beta_0}$ are integrable at $|y|=0$ even for $\delta >0$ close to zero. Furthermore for small $\nu'$ and at a point $|y|>\sqrt{\nu'}$ close to $\sqrt{\nu'}$ we have $|\sqrt{\nu'}|^{1-\epsilon}$ for small $\epsilon>0$ such that the last factor in (\ref{firstderg}) becomes $\exp\left(\frac{-|\nu'|^{1-\epsilon/2}}{8\nu' \sigma} \right)=\exp\left(\frac{-1}{8(\nu')^{\epsilon/2} \sigma} \right)\downarrow 0$ as $\nu'\downarrow 0$. For $\epsilon =2\delta$ we observe that the factor $|\nu'|^{-\delta}$ is damped by this exponential factor such that the upper bound in (\ref{uppb}) is independent of $\nu'$ for these terms.  

More explicitly, for $k=2$ (\ref{Navlerayscheme}) we have
\begin{equation}\label{Navlerayscheme2}
\begin{array}{ll}
 v^{\nu',-,2}_i=v^{\nu',-,0}_i(s_0,.)\ast_{sp}G_{\nu'}
+\sum_{j=1}^D v^{\nu',-,1}_j\frac{\partial v^{\nu',-,1}_i}{\partial x_j}\ast G_{\nu'}\\
\\-\sum_{j,m=1}^D\int_{{\mathbb R}^D}\left( \frac{\partial}{\partial x_i}K_D(.-y)\right) \sum_{j,m=1}^D\left( \frac{\partial v^{\nu',-,1}_m}{\partial x_j}\frac{\partial v^{\nu',-,1}_j}{\partial x_m}\right) (t,y)dy\ast G_{\nu'}\\
\\
=:v^{\nu',-,0}_i(s_0,.)\ast_{sp}G_{\nu'}
+B_1\ast G_{\nu'}-L_1\ast G_{\nu'},
\end{array}
\end{equation}
where $B_1$ and $L_1$ denote abbreviations of the next order of approximation of the Burgers term and the Leray projection term.
The estimate
\begin{equation}
{\big |}v^{\nu',-,0}_i(s_0,.)\ast_{sp} G_{\nu',j}{\big |}\leq cr^{\beta_0-1}.
\end{equation}
Hence
\begin{equation}
{\Big |}\left( v^{\nu',-,1}_j\frac{\partial v^{\nu',-,1}_i}{\partial x_j}\right)(s,.){\Big |}\leq cr^{2\beta_0-1}
\end{equation}
and as $\frac{\partial}{\partial x_i}K_{D,i}(.-y)\sim \frac{1}{r^2}$ for $D=3$ we have
\begin{equation}
{\Big |} \int_{{\mathbb R}^D}K_{D,i}(.-y) \sum_{j,m=1}^D\left( \frac{\partial v^{\nu',-,1}_m}{\partial x_j}\frac{\partial v^{\nu',-,1}_j}{\partial x_m}\right)  (\sigma,y)dy{\Big |}\leq cr^{2(\beta_0-1)+1}.
\end{equation}
Using the $\nu$-independent Gaussian estimates above with $\Delta r:=|x-y|$ 
for the Burgers term $B^1$ we get 
\begin{equation}
{\Big |}B^1\ast D^{\gamma}_xG_{\nu'}(\sigma,.){\Big |}\leq cr^{2\beta_0-1-|\gamma|},
\end{equation}
and  for the Leray projection term $L^0$ we have
\begin{equation}
{\Big |}L^1\ast D^{\gamma}_xG_{\nu'}(\sigma,.){\Big |}\leq cr^{2(\beta_0-1)+1-|\gamma|}.
\end{equation}
Hence, for first order multivariate spatial derivatives $0\leq |\gamma |\leq 1$ and $|x|\geq 1$, for some $\Delta >0$ for all $s \in [s_0,s_0+\Delta]$ there exists a finite constant $c$ which is independent of $\nu >0$ such that for all $\epsilon >0$ small enough we have
\begin{equation}\label{polybound0}
~{\Big |}D^{\gamma}_x \delta v^{\mbox{init},\nu',-,2}_i(\sigma,x){\Big |}\leq  
c r^{2\beta_0-1-\epsilon-|\gamma|},
\end{equation} 
and 
\begin{equation}\label{polybound20}
\forall~ k\geq 3~{\Big |}D^{\gamma}_x \delta v^{\nu',-,k}_i(\sigma,x){\Big |}\leq  cr^{2\beta_0-1-\epsilon-|\gamma|},
\end{equation}
where the constant $c$ is independent of $\nu'$.
Local time contraction results are obtained as for the original Navier Stokes equation, and on some time interval $[s_0,s_0+\Delta]$ we may use representations of solutions of the form
\begin{equation}
\begin{array}{ll}
v^{\nu',-}_i=
v^{\nu',-,0}_i(s_0,.)\ast_{sp}G_{\nu'}+\delta v^{\mbox{init},\nu',-,2}_i\\
\\
+\sum_{k\geq 3}\delta v^{\nu',-,k}_i.
\end{array}
\end{equation}
Hence, for $s\in [s_0,s_0+\Delta]$
\begin{equation}
v^{\nu',-}_i(s,.)\in {\cal C}^{m(D+1)}_{pol,1}
\end{equation}
Using this strong polynomial decay and the compactness argument of the preceding section
we get a subsequence $\nu_k\downarrow 0$ such that
\begin{equation}
e^-\in C^{1,1}~\mbox{where }¸\forall s\in [s_0,s_0+\Delta ]~~e^{-}_i(s,.)=\lim_{k\uparrow \infty}v^{\nu_k,-}_i(s,.)\in {\cal C}^{m(D+1)}_{pol,1}.
\end{equation}

\item[iii)] 
In order to strengthen the regularity result we use the semi-group property of the Euler-and Navier stakes equation operator and the estimates of the preceding item, where we consider the local representation for $\sigma\in [s_0,s_0+\Delta]$
\begin{equation}\label{solrep1}
\begin{array}{ll}
v^{\nu',-}_i(\sigma,.)=
v^{\nu',-,0}_i(s_0,.)\ast_{sp}G_{\nu'}+\delta v^{\mbox{init},\nu',-,2}_i(\sigma,.)\\
\\
+\sum_{k\geq 3}\delta v^{\nu',-,k}_i(\sigma,.).
\end{array}
\end{equation}
For $\sigma >s_0$ the first term in (\ref{solrep1}) $v^{\nu',-,0}_i(s_0,.)\ast_{sp}G_{\nu'}(\sigma,.)$ is a smooth function which implies that for some finite constant $c'$ and  
\begin{equation}\label{polybound21}
\forall~ k\geq 2~{\Big |}D^{\gamma}_x  v^{\nu',-,k}_i(\sigma,x){\Big |}\leq  c'r^{2\beta_0-1-\epsilon-|\gamma|},
\end{equation}
Iterating the argument of the preceding item with initial data $v^{\nu',-,0}_i(s_1,.)$ at $s_1 >s_0$ once we get  
\begin{equation}
{\Big |}L\ast D^{\gamma}_xG_{\nu'}(\sigma,.){\Big |}\leq cr^{4(\beta_0-1)+1-|\gamma|},
\end{equation}
where $L$ denotes the Leray projection operator applied to the local time solution in $[s_1,s_1+\Delta]$ for $\Delta >0$ as above. Similar for the burgers term (where slightly stronger regularity can be proved after one iteration of the regularity argument of item ii) with data $v^{\nu',-,0}_i(s_1,.)$.
We get
\begin{equation}
\forall s\in [s_0,s_0+\Delta ]e^{-}_i(\sigma,.)=\lim_{\nu\downarrow 0}v^{\nu',-}_i(s,.)\in {\cal C}^{m(D+1)}_{pol,m}
\end{equation}
for $s_0+\Delta\geq s_1>s_0$ .

\item[iv)] Choose a time horizon $T>0$. Let $s_0=T$ and consider the time-reversed incompressible Euler equation. From the previous step we have a local solution $e^{\nu,-}_i,~1\leq i\leq D=3$ with data in $H^2$ which correspond to a vorticity blow up at $s_0=T$. We have $e^{\nu,-}_i\in C^{1,1}\left(\left( s_0,s_0+\Delta \right] \right) $ for $1\leq i\leq D$ and for some $\Delta >0$. Moreover as in the previous step such that contraction holds for the higher order increments $\delta v^{\nu,-,k}_i$ with $k\geq 3$ as in (\ref{vel}). Moreover, $e^{-}_i(s,.)\in C^m\cap {\cal C}^{m(D+1)}_{pol,m}$ for $s_0<s\leq s_0+\Delta$. Hence the global solution branch technique of the preceding section can be applied for the time-reversed Euler equation with data $e^{-}_i(s,.)\in C^m\cap {\cal C}^{m(D+1)}_{pol,m}$ for $s>s_0$. It follows that there is a global solution branch $e^{-}_i,~1\leq i\leq D$ defined on the time interval $(s_0,s_0+T]$. Then the time transformation $t=-s_0+T$ implies that $t\rightarrow e_i(t,.)=e^{-}_i,~1\leq i\leq D$ is a global regular solution branch of the incompressible Euler equation  on the time interval $[0,T)$ with a vorticity blow up at time $T$, where $T>0$ is arbitrarily large. 

\end{itemize}

\end{document}